\documentclass[11pt]{article}
\usepackage{epigamath}

\setpapertype{A4}

\usepackage[english]{babel}

\usepackage[dvips]{graphicx}     

\title{On complete reducibility in characteristic $p$}
\titlemark{On complete reducibility in characteristic $p$}
\author{V. Balaji, P. Deligne, and A.\,J. Parameswaran}
\authoraddresses{
\authordata{V. Balaji}{\firstname{V.} \lastname{Balaji}\\
\institution{Chennai Mathematical Institute SIPCOT IT Park, Siruseri-603103, India}\\
\email{balaji@cmi.ac.in}}

\authordata{P. Deligne}{\firstname{P.} \lastname{Deligne}\\
\institution{School of Mathematics, Institute for Advanced Study, Princeton,
NJ 08540, USA}\\
\email{deligne@math.ias.edu}}

\authordata{A.\,J. Parameswaran}{\firstname{A.\,J.} \lastname{Parameswaran}\\
\institution{Tata Institute of Fundamental Research, Mumbai-400005, India}\\
\email{param@math.tifr.res.in}}

\authordata{Z. Yun}{\firstname{Z.} \lastname{Yun}\\
\institution{Department of Mathematics, Stanford University, Stanford, CA 94305, USA}\\
\email{zwyun@math.stanford.edu}}
}
\authormark{V. Balaji, P. Deligne, \&\ A.\,J. Parameswaran}
\date{\vspace{-5ex}} 
\journal{\'Epijournal de G\'eom\'etrie Alg\'ebrique} 
\acceptation{Received by the Editors on November 9, 2016, and in final form on April 12, 2017. \\ Accepted on May 9, 2017.}

\newcommand{\articlereference}{Volume 1 (2017), Article Nr. 3}

\acknowledgement{The research of the first author was partially supported by the J.C. Bose Research grant.}

\usepackage{stmaryrd}
\usepackage{amscd}
\usepackage[all]{xypic}

\usepackage{frcursive}

\newtheorem{guess}{\bf Theorem}[section]
\newcommand{\bth}{\begin{guess}$\!\!\!${\bf }~}
\newcommand{\eeth}{\end{guess}}

\newtheorem{propo}[guess]{\bf Proposition}

\newcommand{\bprop}{\begin{propo}$\!\!\!${\bf }~}
\newcommand{\eprop}{\end{propo}}

\newtheorem{lema}[guess]{\bf Lemma}
\newcommand{\blem}{\begin{lema}$\!\!\!${\bf }~}
\newcommand{\elem}{\end{lema}}

\newtheorem{defe}[guess]{\bf Definition}
\newcommand{\bdefe}{\begin{defe}$\!\!\!${\it }~\rm}
\newcommand{\edefe}{\end{defe}}

\newtheorem{coro}[guess]{\bf Corollary}
\newcommand{\bcor}{\begin{coro}$\!\!\!${\bf }~}
\newcommand{\ecor}{\end{coro}}

\newtheorem{rema}[guess]{\bf Remark}
\newcommand{\brem}{\begin{rema}$\!\!\!${\it }~\rm}
\newcommand{\erem}{\end{rema}}

\newtheorem{assump}[guess]{\bf Assumption}

\newcommand{\spec}{{\rm Spec}\,}
\newtheorem{notation}{Notation}[section]
\newcommand{\bnot}{\begin{notation}$\!\!\!${\bf }~~\rm}
\newcommand{\enot}{\end{notation}}

\newcommand{\bpr}{\begin{proof}}
\newcommand{\epr}{\hfill $\Box$ \end{proof}}

\numberwithin{equation}{guess}

\newcommand{\beqa}{\begin{eqnarray}}
\newcommand{\eeqa}{\end{eqnarray}}

\newtheorem{example}[guess]{\bf Example}
\newtheorem{say}[guess]{\!\!}
\newcommand{\bsem}{\begin{say}\rm}
\newcommand{\esem}{\end{say}}

\newsymbol \bulletledarrowleft 1309

\newcommand{\bz}{{\mathbb Z}}
\newcommand{\bq}{{\mathbb Q}}

\newcommand{\bg}{{\mathbb G}}

\newcommand{\br}{{\mathbb R}}
\newcommand{{\bh}}{{\mathbb H}}

\newcommand{\co}{{\mathcal O}}

\newtheorem{ack}{Acknowledgments\!\!}

\DeclareFontFamily{OT1}{rsfs}{}
\DeclareFontShape{OT1}{rsfs}{n}{it}{<-> rsfs10}{}
\DeclareMathAlphabet{\mathscr}{OT1}{rsfs}{n}{it}

\newcommand{\eL}{{\mathscr L}}

\newcommand{\parallelslant}{\sslash}

\begin{document}


\maketitle

\contribution{with an appendix by Zhiwei Yun}

\dedication{in memoriam Vikram Mehta}

\begin{prelims}

\vspace{-0.1cm}

\def\abstractname{Abstract}
\abstract{Let $G$ be a reductive group over a field $k$ which is algebraically closed
of characteristic $p \neq 0$. We prove a structure theorem for a class of
subgroup schemes of $G$ when $p$ is bounded below by the Coxeter number of $G$. As
applications, we derive semi-simplicity results generalizing earlier results
of Serre proven in 1998 and also obtain an analogue of Luna's \'etale slice
theorem for suitable bounds on $p$.}



\vspace{0.05cm}

\languagesection{Fran\c{c}ais}{%

\textbf{Titre. Compl\`ete r\'eductibilit\'e en caract\'eristique $p$}
\commentskip
\textbf{R\'esum\'e.}
Soit $G$ un groupe r\'eductif sur un corps $k$ qui est alg\'ebriquement clos de ca\-ract\'eristique $p\neq 0$. Nous prouvons un th\'eor\`eme de structure pour une classe de sch\'emas en sous-groupes de $G$ lorsque $p$ est minor\'e par le nombre de Coxeter de $G$. En guise d'applications, nous en d\'eduisons des r\'esultats de semi-simplicit\'e qui g\'en\'eralisent des r\'esultats ant\'erieurs de Serre \'etablis en 1998 et nous obtenons \'egalement un analogue du th\'eor\`eme de Luna (\emph{\'etale slice}) pour des bornes convenables sur $p$.}

\end{prelims}


\newpage

\tableofcontents

\section{Introduction}

We work over an algebraically closed field  $k$ of characteristic $p >0$. 

In \cite{serre2}, Serre showed that if semi-simple representations $V_{_i}$ of a group $\Gamma$ are such that $\sum({\dim}V_{_i} - 1) < p$, then their tensor product is semi-simple. In \cite{mour}, he more generally considers the case where $\Gamma$ is a subgroup of $G(k)$, for $G$ a reductive group, and where $\Gamma$ is $G$-cr, meaning that whenever $\Gamma$ is contained in a parabolic subgroup $P$, it is already contained in a Levi subgroup of $P$. For $G = \prod {\text{GL}}(V_{_i})$, this is equivalent to the semi-simplicity of the representations $V_{_i}$ of $\Gamma$. For  a  representation $V$ of $G$, one can then ask under what conditions does $V$ becomes semi-simple, when considered as a representation of $\Gamma$. In \cite[Theorem 6, page 25]{mour}, Serre shows this is the case when the Dynkin height $ht_{_G}(V)$ is less than  $p$. For $G = \prod {\text{GL}}(V_{_i})$ and $V = \bigotimes V_{_i}$, one has $ht_{_G}(V) = \sum ({\dim}V_{_i} - 1)$.

In \cite{deligne}, the results of \cite{serre2} were generalized to the case when the $V_{_i}$ are semi-simple representations of a group scheme $\mathfrak G$. In this paper, we consider the case when $\mathfrak G$ is a subgroup scheme of a reductive group $G$ and generalize \cite{mour} (see \ref{4.8}) and \cite{deligne} (see \ref{strong}). As in \cite{deligne}, we first have to prove a structure theorem  (\ref{strong}) on doubly saturated  (see \ref{saturated and infsat}) subgroup schemes $\mathfrak G$   of reductive groups $G$. The proof makes crucial use of a result of Zhiwei Yun on root systems. The appendix contains the result.

In Section 5, we consider a reductive group acting on an affine variety $X$ and a point $x$ of $X$ whose orbit $G.x$ is closed in $X$. We prove a schematic analogue of \cite[Propositions 7.4, 7.6]{bardsley} under some conditions on the characteristic of $k$. More precisely, if $X$ embeds in a $G$-module $V$ of low height, then we  obtain, as a consequence of \ref{strong}(2), an analogue of Luna's \'etale slice theorem (\ref{luna}, \ref{5.6}). In \cite{bardsley} the language of schemes was not used and as a consequence the orbit $G.x$ had to be assumed ``separable". An orbit $G.x$ is separable  if and only if the stabilizer $G_{_x}$ is reduced.

\section{Saturation and Infinitesimal saturation}

\bsem\label{nilpotentindex} Let $G$ be a  reductive algebraic group over $k$. Our terminology is that of \cite{sga3}: {\em reductive} implies  {\em smooth and connected}. By an {\em algebraic group} we will mean an affine group  scheme of finite type over $k$.  Fix a maximal torus $T_{_G}$ and a Borel subgroup containing $T_{_G}$ which determines a root system $R$ and a set of positive roots $R^+$.  Let  $\langle~,~ \rangle$ be the natural pairing between the characters and co-characters and for each root $\alpha$, let $\alpha^{^\vee}$ be the corresponding coroot.  

If the root system $R$ associated to $G$ is irreducible, the Coxeter number ${\tt h}_{_G}$ of $R$ and of $G$, admits the following equivalent descriptions:
\begin{enumerate}[\rm (1)]
\item  It is the order of the Coxeter elements of the Weyl group $W$. This shows that $R$ and the dual root system $R^{\vee}$ have the same Coxeter number. 
\item  Let $\alpha_{_0}$ be the highest root and $\sum n_{_i} \alpha_{_i}$ its expression as a linear combination of the simple roots. One has:
\beqa
{\tt h}_{_G} = 1 + \sum n_{_i}.
\eeqa

\item  Applying this to the dual root system, one gets ${\tt h}_{_G} = \langle \rho, \beta^{^\vee} \rangle + 1$ where $\rho$ is half the sum of positive roots and $\beta^{^\vee}$ is the highest coroot. Indeed, $\rho$ is also the sum of the fundamental weights $\omega_{_i}$, and $\langle \omega_{_i}, \alpha_{_j}^{^\vee} \rangle = \delta_{_{ij}}$.

\end{enumerate}
For a general reductive group $G$, define the Coxeter number ${\tt h}_{_G}$ to be the largest among the ones for the irreducible components of $R$. 

It follows from (2) above that if $G$ is a reductive group, $U$ the unipotent radical of a Borel subgroup, and $\mathfrak u:= {\text{Lie}}(U)$, the descending central series of $\mathfrak u$, defined by $Z^{^1} \mathfrak u = \mathfrak u$ and   $Z^{^i} \mathfrak u = [\mathfrak u, Z^{^{i-1}}\mathfrak u]$,  satisfies
\beqa
Z^{^{{\tt h}_{_G}}}(\mathfrak u) = 0.
\eeqa
For the group $U$, similarly $Z^{^{{\tt h}_{_G}}}(U) = (1)$.

The Lie algebra $\mathfrak g$ of $G$ is a $p$-Lie algebra. Let $X \mapsto X^{^{[p]}}$ denote the $p$-power operation on $\mathfrak g$. If  
${\text{\cursive n}}$ in $\mathfrak g$ is nilpotent and if $p \geq {\tt h}_{_G}$, then ${\text{\cursive n}} ^{^{[p]}} = 0$. To check this, we may assume that ${\text{\cursive n}}$ is in the Lie algebra $\mathfrak u$ of the unipotent radical $U$ of a Borel subgroup. For each positive root $\alpha$, let $X_{_\alpha}$ be a basis for  the root subspace $\mathfrak g_{_{\alpha}}$.  Express ${\text{\cursive n}}$ as $\sum_{_{\alpha \in R^{^{+}}}} a_{_{\alpha}}.X_{_{\alpha}}$. Each $X_{_\alpha}$ is the infinitesimal generator of an additive group. It follows that 
$X_{_{\alpha}}^{^{[p]}} = 0$ for each $\alpha$. Observe that
\beqa
{\text{\cursive n}} ^{^{[p]}} =  \sum_{_{\alpha \in R^{^{+}}}} a_{_{\alpha}}^{^{[p]}}.X_{_{\alpha}}^{^{[p]}} ~~~(\mathrm{modulo}~~ Z^{^{p}}\mathfrak u)
\eeqa
and $Z^{^{p}}\mathfrak u = 0$, as $p \geq {\tt h}_{_{G}}$
(see \cite[page 10]{mcninch}).  \esem

\bsem\label{springerstuff} Let  ${\mathfrak g}_{_{nilp}}$ (resp. $G^{^u}$) be the reduced subscheme of   ${\text{Lie}}(G)$ (resp. $G$) with points the nilpotent (resp. unipotent) elements.  Let $U$ be the unipotent radical of a Borel subgroup. For $p \geq {\tt h}_{_{G}}$,  the Campbell-Hausdorff group law $\circ$ makes sense in characteristic $p$ and turns $\mathfrak u := {\text{Lie}}(U)$ into an algebraic group over $k$. This is so since $Z^{^{p}}\mathfrak u = 0$.  Further, there is an unique isomorphism
\beqa\label{chexp}
{\text {exp}}: ({\text{Lie}}(U), \circ) \stackrel{{\sim}}\longrightarrow U
\eeqa
equivariant for the action of $B$ and whose differential at the origin is the identity. If in addition the simply connected covering of the derived group of $G$ is an \'etale covering (which is the case for $p > {\tt h}_{_{G}}$, and could fail when $p = {\tt h}_{_{G}}$ due to the presence of $\text{SL}(p)$ as factors in the covering), then there is a unique $G$-equivariant isomorphism  (\cite[Theorem 3, page 21]{mour} see Section \hyperref[aip]{6} for details):
\beqa\label{springexp}
{\text {exp}}:{\mathfrak g}_{_{nilp}} \to G^{^u}
\eeqa 
which induces \eqref{chexp} on each unipotent radical of a Borel subgroup. Let $\log: G^{^u} \to {\mathfrak g}_{_{nilp}}$ denote its inverse.

For $u$ a unipotent element of $G(k)$, one defines the ``$t$-power map"~$t \mapsto u^t$, from $\bg_{_a}$ to $G$, by  
\beqa\label{springnori}
t \mapsto {\text {exp}}(t~{\text {log}} ~u).
\eeqa

For $G = {\text{GL}}(V)$ such a map $t \mapsto u^t$ is more generally defined for any $u$ in $G$ such that $u^{^{[p]}} = 1$. It is  given by the truncated binomial expression \cite[4.1.1, page 524]{serre2}:
\beqa\label{nori}
t \mapsto u^{t} := \sum_{_{i < p}} \binom{t}{i}(u-1)^{^i}.
\eeqa
Similarly, $X$ in ${\text{End}}(V)$ such  that $X^{^{[p]}} = 0$ defines a morphism $t \mapsto \exp(tX)$ from $\bg_{_a}$ to ${\text{GL}}(V)$ given by the truncated exponential series:
\beqa\label{exp}
t \mapsto  I + tX + \frac{(tX)^{^2}}{2 !} + \ldots + \frac{(tX)^{^{p-1}}}{{(p-1)} !}.
\eeqa
\esem

\noindent
{\em  Until the end of \S3, we make the following assumptions  on the reductive group $G$}.

\begin{assump}\label{basicassumptions} \rm { Let  $\widetilde{G}$ be the simply connected covering of the derived group $G'$. We assume that $p \geq {\tt h}_{_{G}}$ and that the map $\widetilde{G} \to G'$ is \'etale}.   \end{assump}

In particular, by \ref{nilpotentindex} and \ref{springerstuff}, the exponential map \eqref{springexp} is defined, every unipotent element in $G(k)$ is of order $p$, and every nilpotent in ${\text{Lie}}(G)$ is a $p$-nilpotent. One can then define the notions of saturation and infinitesimal saturation  of subgroup schemes $\mathfrak G \subset G$ as follows (see Remark \ref{thepgl(p)case} for the case when $G = {\text{PGL}}(p)$). 

\bdefe\label{saturated and infsat}(\cite[\S4]{serre2}, \cite[Definition 1.5]{deligne}) 
\begin{enumerate}[\rm (1)]
\item A subgroup scheme $\mathfrak G \subset G$ is called {\em saturated} if  for every $u$ in $\mathfrak G(k)$ which is unipotent, the homomorphism $t \mapsto u^{t}$ \eqref{springnori} from $\bg_{_a}$ to $G$ factors through $\mathfrak G$. 
\item A subgroup scheme $\mathfrak G \subset G$ is called {\em infinitesimally saturated} if for every nilpotent $X$ in ${\text{Lie}}(\mathfrak G)$, the morphism $t \mapsto \exp(tX)$ \eqref{springexp} from ~$~\bg_{_a}$ to $G$ factors through $\mathfrak G$. 

\item $\mathfrak G$ is {\em doubly saturated} if it is saturated and infinitesimally saturated.
\end{enumerate}
\edefe

An element of ${\text{Lie}}(\mathfrak G)$ is nilpotent if and only if it is nilpotent as an element of ${\text{Lie}}(G)$. The reference to the exponential map \eqref{springexp} in (2) therefore makes sense. One way to see that the notions
  of nilpotence for elements of ${\text{Lie}}(G)$ and ${\text{Lie}}(\mathfrak{G})$ coincide
  is to observe that the inclusion of ${\text{Lie}}(\mathfrak G)$ in ${\text{Lie}}(G)$ is a morphism of  $p$-Lie algebras and that $X$ is nilpotent if and only if it is killed by an iterated $p$-power map, i.e. $X^{^{[p^\ell]}} = 0$. 

Let $\mathfrak G^{^0}$ be the identity component of $\mathfrak G$ and $\mathfrak G^{^0}_{_{red}}$ the reduced subscheme of $\mathfrak G^{^0}$.

\bth\label{strong} Let ${\mathfrak G} \subset G$ be a $k$-subgroup scheme which  is infinitesimally saturated. Assume that  if $p~=~{\tt h}_{_G}$,  $\mathfrak G^{^0}_{_{red}}$ is reductive.  Then 
\begin{enumerate}[\rm (1)]
\item The group  $\mathfrak G^{^0}_{_{red}}$ and  its unipotent radical $R_{_u}(\mathfrak G^{^0}_{_{red}})$ are normal subgroup schemes of $\mathfrak G$ and the quotient group scheme $\mathfrak G^{^0}/\mathfrak G^{^0}_{_{red}}$ is of multiplicative type.
\item If $\mathfrak G^{^0}_{_{red}}$ is reductive, there exists a central, connected subgroup scheme of multiplicative type $M \subset \mathfrak G^{^0}$ such that the morphism $M \times \mathfrak G^{^0}_{_{red}} \to \mathfrak G^{^0}$ realizes $\mathfrak G^{^0}$ as a quotient of $M \times \mathfrak G^{^0}_{_{red}}$. 
\end{enumerate}
\eeth

We note that since $k$ is assumed to be algebraically closed, a group scheme of multiplicative type is simply a {\em diagonalisable} group scheme.

The proof of part (2) of \ref{strong} will occupy most  section 2, until \ref{lieseq}. Part (1) will be proven in section 3.

By \cite[Lemma 2.3]{deligne}, the conclusions of \ref{strong} hold for $\mathfrak G$ if and only if they hold for the identity component $\mathfrak G^{^o}$. { Until the end of section 3, we will assume that $\mathfrak G$ is connected}.

\blem\label{2.6} If $\mathfrak G$ is an infinitesimally saturated subgroup scheme of $G$, every nilpotent element ${\text{\cursive n}}$ of $\text{Lie}(\mathfrak G)$ is in $\text{Lie}(\mathfrak G_{_{red}})$. \elem

\bpr As $\bg_{_a}$ is reduced, the morphism$t \mapsto \exp(t {\text{\cursive n}})$ maps $\bg_{_a}$ to $\mathfrak G_{_{red}} \subset G$. Identifying $\text{Lie}(\bg_{_a})$ with a copy of the ground field, the
   image of $1$ in $\text{Lie}(\bg_{_a})$ is ${\text{\cursive n}}$. \epr

Part (2) of \ref{strong} is a corollary of \ref{2.6} and of the following theorem, which does not refer to $G$ anymore. 
 
\bth\label{2.7} Let $\mathfrak G$ be a connected algebraic group such that
\begin{enumerate}[\rm (a)]
\item $\mathfrak G_{_{red}}$ is reductive,
\item any nilpotent element of $\text{Lie}(\mathfrak G)$ is in $\text{Lie}(\mathfrak G_{_{red}})$.
\end{enumerate}
Then, the conclusion of \ref{strong}(2) holds. As a consequence, $\mathfrak G_{_{red}}$ is a normal subgroup scheme of $\mathfrak G$ and $\mathfrak G/\mathfrak G_{_{red}}$ is of multiplicative type. \eeth

 Let $T$ be a maximal torus of $\mathfrak G_{_{red}}$, and let $H$ be the centralizer of $T$ in $\mathfrak G$. One has $H \cap \mathfrak G_{_{red}} = T$. It follows from (b) that any nilpotent element ${\text{\cursive n}}$ of $\text{Lie}(H)$ is in $\text{Lie}(T)$, hence vanishes. By the following lemma, $H$ is of multiplicative type, and in particular commutative. 

\blem Let $H$ be a connected algebraic group over $k$. If all the elements of ${\text{Lie}}(H)$ are semi-simple, then $H$ is of multiplicative type.\elem

\bpr (See also \cite[IV, \S3, Lemma 3.7]{demgab}.)
\underline{${\text{Lie}}(H)$ is commutative}: Fix $x$ in ${\text{Lie}}(H)$, and let us show that it is central in ${\text{Lie}}(H)$. As $ad x$ is semi-simple, it suffices to show that if $y$ is in an eigenspace of $ad x$, i.e. $[x,y] = \lambda y$, then $x$ and $y$ commute, i.e. $\lambda = 0$. Let $W$ be the vector subspace of ${\text{Lie}}(H)$ generated by the $y^{^{[p^\ell]}}$ $(\ell \geq 0)$. The $y^{^{[p^\ell]}}$ commute. The map $z \mapsto z^{^{[p]}}$ therefore induces a $p$-linear map from $W$ to itself, injective by assumption. This implies that $W$ has a basis $e_{_i}, (1\leq i \leq N)$ consisting of elements such that $e_{_i}^{^{[p]}} = e_{_i}$, and $(\sum a_{_i} e_{_i})^{^{[p]}} = \sum (a_{_i})^{^{p}} e_{_i}.$ To see this, we view $W$ as an algebraic group. The morphism $W \to W$ given by $x \mapsto x^{^{[p]}} - x$ is \'etale. Let $K$ be its kernel. In any basis it is defined by $d = \text{dim}(W)$ equations of degree $p$, namely: $(x^{^{[p]}})_{_i} - x_{_i} = 0, i = 1, 2, \ldots, d$, with an \'etale set of solutions and no solution at $\infty$ (the system of homogeneous equations $(x^{^{[p]}})_{_i} = 0$ has no non-zero solutions by assumptions). By Bezout, the group of solutions is hence a $(\bz/p)^{^d}$. Let $e_{_i}$ be a basis for it. If the $e_{_i}$ were not linearly independent in $W$, there would be a linear dependence $\sum _{_{i \in J}} a_{_i} e_{_i} = 0$, with the $a_{_i} \in k$ non-zero, involving a minimal number of $e_{_i}$. We also have $\sum a_{_i}^{^p} e_{_i} = 0$. Hence, by the minimality of the $e_{_i}$, for some $\lambda$ and $\forall i$, $a_{_i}^{^p} = \lambda a_{_i}$. Rescaling by $\lambda^{^{\frac{1}{p-1}}}$, we get $a_{_i} \in \mathbb F_{_p}$, which is a contradiction.

\blem For $b = (b_{_i})$ in $k^N$, define $b^{^{[p]}}:= (b_{_i}^{^{[p]}})$. Then, any $b$ in $k^N$ is a linear combination of the $b^{^{[p^\ell]}}$ for $\ell \geq 1$.\elem

\bpr The $b^{^{[p^a]}}$ $(a \geq 0)$ are linearly dependent. A linear dependence relation can be written  
\beqa
\sum_{_{j \geq m}} c_{_j} b^{^{[p^j]}}= 0
\eeqa
with $c_{_m} \ne 0$. Extracting $p^{^m}$-roots, we get  
\beqa
(\sum_{_{j \geq 0}} d_{_j} b^{^{[p^j]}}) = 0
\eeqa
where $d_{_j} = c_{_{m+j}}^{1/p^m}$. In particular, $d_{_0} \ne 0$, proving the lemma.
\epr

\noindent
\underline{End of proof of commutativity}: From the lemma above, $y$ is a linear combination of the $y^{^{[p^\ell]}}$ $(\ell > 0)$. The bracket $[y^{^{[p^\ell]}}, x]$ vanishes for $(\ell > 0)$. Indeed, it is $(ad y)^{^{[p^\ell]}}(x)$ which vanishes because $[y,[y,x]] = [y, -\lambda y] = 0$. It follows that $[y,x]$ vanishes too.

The $p$-Lie algebra ${\text{Lie}}(H)$ is commutative and hence its $p$-power operation ${\text{Lie}}(H) \to {\text{Lie}}(H)$ is injective. It therefore has a basis $e_{_i}$ such that $e_{_i}^{^{[p]}} = e_{_i}$. The dual of its restricted enveloping algebra is the affine bigebra of  the kernel $K$  of the Frobenius morphism $F:H \to H^{^{(p)}}$, where $H^{^{(p)}}$ is obtained from $H$ by extension of scalars $\lambda \mapsto \lambda^{^{[p]}}$, $k \to k$. It follows that $K$ is a product of $\mu_{_p}$'s.

The same holds for $H^{^{(p)}}$, which is obtained from $H$ using an automorphism of $k$. The same holds for each $H^{^{(p^\ell)}}$ as well.

For any $n$, the kernel $K_{_n}$ of the iterated Frobenius map $F^{^n}: H \to H^{^{(p^\ell)}}$ is an iterated extension of subgroups of the kernel of the Frobenius of the $H^{^{(p^i)}}$, $i < n$. It is hence of multiplicative type, being an iterated extension of connected groups of multiplicative type \cite[XVII, 7.1.1]{sga3}.

The $K_{_n}$ form an increasing sequence. By the proof of \cite[Proposition 1.1]{deligne}, there exist subgroup of multiplicative type $M$ of $H$ containing all $K_{_n}$. As $H$ is connected and as $M$ contains all infinitesimal neighbourhoods of the identity element, one has $M = H$.
\epr

Whenever a group $M$ of multiplicative type acts on a group $K$, its action on ${\text {Lie}}(K)$ defines a weight decomposition:
\beqa\label{weightdecomp}
{\text {Lie}}(K) = \bigoplus_{_{\beta \in X(M)}} {\text {Lie}}(K)^{^{\beta}}.
\eeqa
If $v \in {\text {Lie}}(K)^{^{\beta}}$, then $v^{^{[p]}}$ is in ${\text {Lie}}(K)^{^{p.\beta}}$. Indeed, after any extension of scalars $R/k$, if $m$ is in $M(R)$, 
\beqa
m(v^{^{[p]}}) = (m(v))^{^{[p]}} = (\beta(m)v)^{^{[p]}} = \beta(m)^{^{[p]}} v^{^{[p]}}.
\eeqa

\blem\label{2.10} Let $M$ be a multiplicative group acting on a group $K$ and let  \eqref{weightdecomp} be the weight space decomposition for the corresponding action on ${\text {Lie}}(K)$. If $\beta \in X(M)$ is not torsion, $\text{Lie}(K)^{^\beta}$ consists of nilpotent elements. \elem

\bpr Indeed, if $\beta$ is not torsion, the $p^{^\ell}\beta$ are all distinct and ${\text {Lie}}(K)^{^{p^{^\ell}\beta}}$ must vanish for $\ell \gg 0$. It follows that the elements of ${\text {Lie}}(K)^{^{\beta}}$ are nilpotent.\epr

Let us apply this to the action of $T$ on $\mathfrak G$ by inner automorphisms.

\blem\label{2.11} If $\beta$ in $X(T)$ is not zero, the weight spaces ${\text {Lie}}(\mathfrak G)^{^\beta}$ equals ${\text {Lie}}(\mathfrak G_{_{red}})^{^\beta}$.\elem

\bpr Indeed, ${\text {Lie}}(\mathfrak G)^{^\beta}$ consists of nilpotent elements. By our assumption \ref{2.7}(b), it is contained in ${\text {Lie}}(\mathfrak G_{_{red}})$. \epr

Let $B$ a Borel subgroup of $\mathfrak G_{_{red}}$ containing $T$, and let $U$ be its unipotent radical.

\blem\label{3.8} Under the assumption of \ref{2.7}, $H$ normalizes $U$. \elem

\bpr Let $C$ in $X(T) \otimes \br$ be the cone generated by the positive roots  relative to $B$ and define
\beqa
C^{*} := C \setminus \{0\}.
\eeqa
As in \ref{2.11}, we let $T$ act on $\mathfrak G$ by conjugation ($t$ acts by $g \mapsto tgt^{^{-1}}$). This action induces actions on $\mathfrak g := {\text {Lie}}(\mathfrak G)$, the affine algebra $A$ of $\mathfrak G$, its augmentation ideal $\mathfrak m$ (defining the unit element), and the dual $\mathfrak g^{\vee}= \mathfrak m/\mathfrak m^{^2} $ of $\mathfrak g$. Similarly $T$ acts on $U$, its Lie algebra $\mathfrak u$, its affine algebra $A_{_U}$ and its augmentation ideal $\mathfrak m_{_U}$. For the action on affine algebras, $t$ in $T$ transforms $f(g)$ into $f(t^{^{-1}}gt)$. From these actions, we get $X(T)$-gradings. By \ref{2.7}(b)
\beqa
\mathfrak u = \bigoplus_{_{\beta \in C^*}}  \mathfrak g^{^{\beta}}.
\eeqa
It follows that $\mathfrak u^{\vee} = \mathfrak m_{_U}/\mathfrak m_{_U}^{^2}$ is the sum of the $(\mathfrak g^{\vee})^{^\beta}$,  for $\beta$ in the negative $-C^*$ of $C^*$. For $n > 0$, the weights by which $T$ acts on $\mathfrak m_{_U}^{^n}/\mathfrak m_{_U}^{^{n+1}}$ are in $-C^*$. As $U$ is connected, the intersection of the $\mathfrak m_{_U}^{^n}$ is reduced to $0$, and the weights by which $T$ acts on $\mathfrak m_{_U}$ are also in $-C^*$.

Let $I$ be the ideal of $A$ generated by the graded components $\mathfrak m^{^\beta}$ of $\mathfrak m$ for $\beta$ not in $-C^*$, and put $A_{_1}:= A/I$. The image in $\mathfrak m_{_U}$ of a $\mathfrak m^{^\beta}$ as above is contained in $\mathfrak m_{_U}^{^\beta}$, hence vanishes. It follows that $U$ is contained in the closed subscheme $U_{_1} = \spec A_{_1}$ of $\mathfrak G$ defined by $I$. As the graded component  defined $\mathfrak m^{^0}$ of $\mathfrak m$ is contained in $I$, the graded component $A_{_1}^{^0}$ of $A_{_1}$ (the $T$-invariants) is reduced to the constants. As $T$ is connected it follows that $U_{_1}$ is connected. As the image of $\mathfrak m_{_U}^{^\beta}$ in $\mathfrak m/\mathfrak m^{^2}$ is $(\mathfrak m/\mathfrak m^{^2})^{^\beta}$, the image of $I$ in $\mathfrak m/\mathfrak m^{^2} = \mathfrak g^{\vee}$ is the orthogonal of $\mathfrak u \subset \mathfrak g$ and the tangent space at the origin of $U_{_1}$ is $\mathfrak u$.

{\em Claim}: The subscheme $U_{_1}$ of $\mathfrak G$ is a subgroup scheme, i.e. the coproduct $\Delta: A \to A \otimes A, f(g) \mapsto f(gh)$ maps $I$ to $I \otimes A \oplus A \otimes I$. Indeed, $\Delta$ respects the gradings; for $k \subset A$ the constants, one has $A = k \oplus \mathfrak m$, and $\Delta$ maps $\mathfrak m$ to $(k \otimes \mathfrak m) \oplus (\mathfrak m \otimes k) \oplus  (\mathfrak m \otimes \mathfrak m)$, and hence $\mathfrak m^{^\beta}$ maps to 
\beqa
(k \otimes \mathfrak m^{^\beta}) \oplus (\mathfrak m^{^\beta}\otimes k) \oplus  \sum_{_{\beta = \beta' + \beta''}} \mathfrak m^{^{\beta'}} \otimes \mathfrak m^{^{\beta''}}.
\eeqa
As $-C^*$ is stable by addition, if $\beta = \beta' + \beta''$ and that $\beta$ is not in $-C^*$, one of $\beta'$ or $\beta''$ is not in $-C^*$, and the corresponding $\mathfrak m^{^{\beta'}}$ or $\mathfrak m^{^{\beta''}}$ is contained in $I$. The claim follows.

To summarize, $U_{_1}$ is connected, and the inclusion $U\subset U_{_1}$ induces an isomorphism ${\text{Lie}}(U) \stackrel{\sim}\longrightarrow {\text{Lie}}(U_{_1})$. As $U$ is smooth, this implies that $U = U_{_1}$.

Since $H$ centralizes $T$, the ideal $I$ is stable by $H$, meaning that $H$ normalizes $U$. 
\epr

\bcor\label{3.9} Under the assumption of \ref{2.7}, $H$ normalizes $\mathfrak G_{_{red}}$.\ecor

\bpr Let $B^{-}$ be the Borel subgroup of $\mathfrak G_{_{red}}$ containing $T$ and opposite to $B$, and let $U^-$ be its unipotent radical. As $U^-$, $T$ and $U$ are normalized by $H$, the big cell $U^- T U \subset \mathfrak G_{_{red}}$ is stable by the conjugation action of $H$, and so is its schematic closure $\mathfrak G_{_{red}}$. \epr

In what follows, we identify schemes with the corresponding fppf sheaves. A quotient such as $\mathfrak G/H$ represents the quotient of the sheaf of groups $\mathfrak G$ by the subsheaf $H$: $\mathfrak G$ is a $H$-torsor over $\mathfrak G/H$.

\blem\label{3.11} The morphism of schemes 
\beqa\label{3.11.1}
\mathfrak G_{_{red}} /T \to \mathfrak G/H
\eeqa
is an isomorphism.\elem

\bpr  As $T$ is the intersection of $\mathfrak G_{_{red}}$ and $H$, the morphism \eqref{3.11.1}, as a morphism of fppf sheaves, is injective. Testing on $\spec(k)$ and $\spec(k[\epsilon]/(\epsilon^{^2}))$, one sees that it is bijective on points and injective on tangent space at each point. It is hence radicial and unramified. Therefore, on some open set of $\mathfrak G_{_{red}} /T$, this morphism is an immersion.  The $\mathfrak G_{_{red}}$ homogeneity then shows that it is a closed embedding. As $\mathfrak G_{_{red}} /T$ is smooth and \eqref{3.11.1} is bijective on $k$-points, to prove that the closed embedding \eqref{3.11.1} is an isomorphism, it suffices to prove that at each point it induces an isomorphism of tangent spaces. By homogeneity, it suffices to check this at the origin. The tangent space at the origin of $\mathfrak G_{_{red}} /T$ is ${\text{Lie}}(\mathfrak G_{_{red}}) /{\text{Lie}}(T)$.

For the tangent space of $\mathfrak G/H$ we proceed as follows. Let  $J$ be the ideal defining $H$ in $\mathfrak G$ and $\eta$ be the ideal defining the origin in $\mathfrak G/H$. As $\mathfrak G$ is a $H$-torsor on $\mathfrak G/H$, the pull-back of $\eta/\eta^{^2}$ to $\mathfrak G/H$ is $J/J^{^2}$, and $\eta/\eta^{^2}$ is the fiber if $J/J^{^2}$ at the origin $e$. By \cite[2.15]{deligne} applied to $T$ acting on $\mathfrak G$ by conjugation, we have an isomorphism $(J/J^{^2})_{_e} \stackrel{\sim}\longrightarrow \bigoplus_{_{\beta \ne 0}} (\mathfrak m/\mathfrak m^{^2})^{^\beta} = ({\text{Lie}}(\mathfrak G_{_{red}}) /{\text{Lie}}(T))^{\vee}$. Taking duals, we see that  \eqref{3.11.1} is an isomorphism near the origin, hence everywhere, proving \ref{3.11}.
\epr

\bcor\label{3.10} Under the assumptions of \ref{2.7}, $\mathfrak G_{_{red}}$ is normal in $\mathfrak G$, and there exists in $\mathfrak G$  a central connected subgroup scheme of multiplicative type $M$ such that the morphism $M \times \mathfrak G_{_{red}} \to \mathfrak G$ realizes $\mathfrak G$ as a quotient of $M \times \mathfrak G_{_{red}}$. \ecor

\bpr By \ref{3.11}, the product map $\mathfrak G_{_{red}} \times H \to \mathfrak G$ is onto, as a morphism of sheaves. As both $\mathfrak G_{_{red}}$ and $H$ normalize $\mathfrak G_{_{red}} \subset \mathfrak G$, so does $\mathfrak G$.

To complete the proof of \ref{3.10} (and thereby of \ref{2.7}), we follow \cite[\S2.25]{deligne}. Let $M$ be the subgroup of $H$ which centralizes  $\mathfrak G_{_{red}}$ (as always, in the scheme theoretic sense). By \cite[Corollary 2.4 (a) , page 476]{demgab}, $M$ is also multiplicative.  Since $\mathfrak G_{_{red}}$ is reductive, the group scheme ${\text {Aut}}_{_{T}}(\mathfrak G_{_{red}})$ of automorphisms which preserve $T$ is  precisely $T^{^{ad}}$, the image of $T$ in the adjoint group. Hence the conjugation action of $H$ on $\mathfrak G_{_{red}}$ gives the exact sequence:
\beqa
1 \to M \to H \to T^{^{ad}} \to 1
\eeqa
and $T$ surjects onto $T^{^{ad}}$ implying that $M$ and $T$ generate $H$. Since $M$ is generated by $M_{_{red}} \subset \mathfrak G_{_{red}}$ and $M^{^0}$, and since $H$ and $\mathfrak G_{_{red}}$ generate $\mathfrak G$ we see that $M^{^0}$ and $\mathfrak G_{_{red}}$ generate $\mathfrak G$. Moreover, $M$ is central. Thus 
\beqa
M^{^0} \times \mathfrak G_{_{red}} \to \mathfrak G
\eeqa
is an epimorphism. This concludes the proof of \ref{2.7} and in particular, \ref{strong}(2). \epr

\blem\label{lieseq}  Suppose  that $H$ is a maximal connected subgroup scheme of multiplicative type of an  algebraic group ${\mathfrak G}$. Let $Z^{^0}_{_{\mathfrak G}}(H)$ be the identity component of the centralizer of $H$ in $\mathfrak G$ and define $U:= Z^{^0}_{_{\mathfrak G}}(H)/H$. Then the sequence:
\beqa\label{basiclieseq}
0 \to {\text {Lie}}(H) \to {\text {Lie}}(Z^{^0}_{_{\mathfrak G}}(H)) \to {\text {Lie}}(U) \to 0
\eeqa
associated to the central extension
\beqa\label{basicseq}
1 \to H \to Z^{^0}_{_{\mathfrak G}}(H) \to U \to 1
\eeqa 
is  exact.
\elem

\bpr Left exactness of \eqref{basiclieseq} is clear. The maximality of $H$  implies that $U$ is unipotent (see \cite[\S2.5, page 590]{deligne}).

Embed $H$ in $\bg_{_m}^r$ as a subgroup scheme. The quotient $H''= \bg_{_m}^r/H$ is a torus, being a quotient of one. The central extension  \eqref{basicseq},
by a push forward, gives a central extension \cite[Expos\'e XVII, Lemma 6.2.4]{sga3}
\beqa\label{basicseq2}
1 \to \bg_{_m}^r  \to E \to U \to 1
\eeqa
and a diagram of groups:
\beqa
\begin{CD}
  \xymatrix{  & 1 \ar[d] & 1 \ar[d] & 1 \ar[d] \\
  1 \ar[r] &H  \ar[r]\ar[d] & Z^{^0}_{_{\mathfrak G}}(H) \ar[r]\ar[d]& U \ar[r]\ar[d]^{=} & 1\\
1 \ar[r] &\bg_{_m}^r  \ar[r] \ar[d]& E \ar[r]\ar[d]  & U \ar[r] \ar[d]_{} & 1\\
1 \ar[r]  & H'' \ar[r]^{=} \ar[d]  & H'' \ar[r]  \ar[d]&1\\
 & 1 & 1 } 
\end{CD}
\eeqa
Since $H''$ is multiplicative  every  nilpotent in ${\text {Lie}}(E)$ maps to $0$ in  ${\text {Lie}}(H'')$ and hence comes from a  nilpotent in ${\text {Lie}}(Z^{^0}_{_{\mathfrak G}}(H))$.

Since $\bg_{_m}^r$ is smooth, by \cite[Expos\'e VII, Proposition 8.2]{sga3} the sequence \eqref{basicseq2} gives an exact sequence 
\beqa\label{basiclieseq2}
0 \to {\text {Lie}}(\bg_{_m}^r) \to {\text {Lie}}(E) \to {\text {Lie}}(U) \to 0.
\eeqa
Since $U$ is unipotent, any element $z$ in ${\text {Lie}}(U)$ is nilpotent. Let $z'$ in ${\text {Lie}}(E)$ be a lift of $z$. The Jordan decomposition makes sense for any $p$-Lie algebra over a perfect field $k$ and uses only the $p$-power map (see for example \cite[Corollary 4.5.9, page 135]{winter2}. Thus, by using the Jordan decomposition of the lift $z'$ in ${\text {Lie}}(E)$ and noting that the semi-simple part gets mapped to zero in ${\text {Lie}}(U)$, we can assume that $z'$ can  also be chosen to be nilpotent. 

Since every  nilpotent in ${\text {Lie}}(E)$ comes from a  nilpotent in ${\text {Lie}}(Z^{^0}_{_{\mathfrak G}}(H))$, we conclude that $z$ gets lifted to a nilpotent in ${\text {Lie}}(Z^{^0}_{_{\mathfrak G}}(H))$. This implies that \eqref{basiclieseq} is also right exact. \epr

\blem\label{infsatofcentralizer} Suppose that ${\mathfrak G}$ is a subgroup scheme of $G$ which  is infinitesimally saturated in $G$. Then the subgroup scheme $Z^{^0}_{_{\mathfrak G}}(H)$ is infinitesimally saturated in $G$; in particular, every non-zero nilpotent in ${\text {Lie}}(Z^{^0}_{_{\mathfrak G}}(H))$ lies in ${\text {Lie}}(Z^{^0}_{_{\mathfrak G}}(H)_{_{red}})$. \elem

\bpr If ${\text{\cursive n}}$ in ${\text {Lie}}(Z^{^0}_{_{\mathfrak G}}(H))$ is nilpotent, the map $\rho:t \mapsto {\text {exp}}(t{\text{\cursive n}}):\bg_{_a} \to G$ factors through $\mathfrak G$. Since $H$ is central in $Z^{^0}_{_{\mathfrak G}}(H)$, the action of $H$ by inner automorphisms fixes ${\text{\cursive n}}$. Since the map ``${\text {exp}}$" is compatible with conjugation  the entire curve $\rho$ is fixed by $H$. Therefore, $\rho$ factors through $Z_{_{\mathfrak G}}(H)$ and hence through $Z^{^0}_{_{\mathfrak G}}(H)_{_{red}}$, since $\bg_{_a}$ is reduced and connected, and hence ${\text{\cursive n}} \in {\text {Lie}}(Z^{^0}_{_{\mathfrak G}}(H)_{_{red}}).$ This completes the proof of the lemma.\epr

We have the following extension of \cite[Lemma 2.7, Corollary 2.11]{deligne}. 

\blem\label{lemma 2.7} Let ${\mathfrak G}$ be as in \ref{infsatofcentralizer}. Let $H \subset \mathfrak G$ be  a maximal connected subgroup scheme of multiplicative type.  Then the central extension \eqref{basicseq}
splits and $U$ is a smooth unipotent group. 
\elem

\bpr  We claim that $U$ is {\em smooth}. By Lemma \ref{lieseq}, every element $z$ in ${\text {Lie}}(U)$ comes from a nilpotent in ${\text{\cursive n}}$ in ${\text {Lie}}(Z^{^0}_{_{\mathfrak G}}(H))$.  By \ref{infsatofcentralizer}, 
${\text{\cursive n}}$  is in ${\text {Lie}}(Z^{^0}_{_{\mathfrak G}}(H)_{_{red}})$. Its image $z$ is hence in ${\text {Lie}}(U_{_{red}})$. Thus, ${\text {Lie}}(U) = {\text {Lie}}(U_{_{red}})$ proving the claim.

Now we have a central extension \eqref{basicseq} with the added feature that $U$ is smooth. By \cite[Expos\'e XVII, Theorem 6.1.1]{sga3} it follows that \eqref{basicseq} splits (uniquely) and  
\beqa\label{splittingtheseq}
Z^{^0}_{_{\mathfrak G}}(H) = H \times U.
\eeqa
\epr

\brem \label{thepgl(p)case} Recall that  in  \ref{springerstuff}, for the existence of an exponential map \eqref{springexp} we assumed that the covering morphism $\widetilde{G} \to G'$ is \'etale, which therefore became a part of the standing assumption (\ref{basicassumptions}). The case which gets excluded is when the simply connected cover of the derived group has factors of ${\text{SL}}(p)$. For instance, when $G'$ is simple with $p = {\tt h}_{_G}$, the only case excluded is $G' = {\text{PGL}}(p)$.

Let $G$ be a reductive group for which $p = {\tt h}_{_G}$. The tables of Coxeter numbers of simple groups show that except for type $A$, where ${\tt h}_{_{\text{SL}(n)}} = n$, Coxeter numbers are even greater than $2$. It follows that $\widetilde{G}$  is a product of ${\text{SL}}(p)$'s and other simple factors with $p$ larger than their Coxeter numbers.  Even when the morphism $\widetilde{G} \to G'$  is {\em not \'etale}, we can still define the notions of {\em saturation} (resp {\em infinitesimal saturation}) of subgroup schemes $\mathfrak G \subset G$ as follows. 

Say $\mathfrak G \subset G$ is {\em saturated} (resp {\em infinitesimally saturated}) if the inverse image of $\mathfrak G$ in $\widetilde{G}$ is saturated (resp infinitesimally saturated). With this definition, Theorem \ref{strong} remains true for $p \geq {\tt h}_{_G}$. 

This notion of saturation (resp infinitesimal saturation) can also be seen in terms of suitably defined $t$-power maps and exponential maps. We restrict ourselves to the case when $G' = {\text{PGL}}(p)$. 
At the level of Lie algebras, the induced morphism:
\beqa
{\mathfrak{sl}}(p) \to {\mathfrak{pgl}}(p)
\eeqa
is a radicial map on the locus  of nilpotent elements. If $A \in {\mathfrak{gl}}(p)$ is a matrix representing an element of ${\mathfrak{pgl}}(p)$, it is nilpotent if all but the constant  coefficients of the characteristic polynomial vanishes, i.e. $Tr(\wedge^{i}(A)) = 0, \forall 1 \leq i < p$ and the characteristic polynomial reduces to $T^{^{p}} - {\text{det}}(A)$. This condition is stable under $A \mapsto A + \lambda.I$ as $(T - \lambda)^{^{p}} - {\text{det}}(A) = T^{^{p}} - ({\text{det}}(A) + \lambda^{^{p}}) = T^{^{p}} - {\text{det}}(A + \lambda.I)$.

We get the unique lift $\widetilde{A} \in  {\mathfrak{sl}}(p)$ by taking $\widetilde{A}:= A - {\text{det}}(A)^{^{1/p}}.I$. Now one has $Tr(\wedge^{i}(\widetilde{A})) = 0$ for $1 \leq i \leq p$ and hence we can define the exponential morphism $\bg_{_a} \to {\text{PGL}}(p)$ as:
\beqa\label{expinpglp}
t \mapsto {\text{exp}}(tA) := {\text{exp}}(t\widetilde{A}).
\eeqa
Likewise, in the case of unipotents in the group ${\text{PGL}}(p)$, the restriction of this map to the locus of unipotent elements:
\beqa
{\text{SL}}(p)^{^u} \to {\text{PGL}}(p)^{^u}
\eeqa
is radicial. Thus, any unipotent $u \in {\text{PGL}}(p)$ has a unique unipotent lift $\widetilde{u} \in {\text{SL}}(p)$ and one can define the $t$-power map $\bg_{_a} \to {\text{PGL}}(p)$ as:
\beqa\label{tpowerinpglp}
t \mapsto u^{^t}:= \mathrm{image}(\widetilde{u}^{^t})
\eeqa
With the notions of $t$-power map and exponential morphisms in place, we can define the notions of {\em saturation} (resp {\em infinitesimal saturation}) of subgroup schemes $\mathfrak G \subset G'$ exactly as in Definition \ref{saturated and infsat} using \eqref{tpowerinpglp} (resp \eqref{expinpglp}) and these coincide with the definitions made above. 

We note however that these ``punctual" maps, i.e. defined for each nilpotent $A$ (resp each unipotent $u$), are not induced by a morphism from ${\mathfrak{pgl}}(p)_{_{nilp}}$ to ${\text{PGL}}(p)^{^u}$ (resp. $\mathbb A^{1} \times {\text{PGL}}(p)^{^u} \to {\text{PGL}}(p)^{^u}: (t,u) \mapsto u^{^t}$ is not a morphism).
\erem

\section{Completion of proof of  the structure theorem \ref{strong}}

We begin by stating a general result on root systems whose proof is given in the appendix.

\bprop\label{yun} Let $R$ be an irreducible root system  with Coxeter number ${\tt h}$ and let $X$ be the lattice spanned by $R$. Let $\phi:X \to \br/\bz$ be a homomorphism. Then there exists a basis $B$ for $R$ such that if $\alpha \in R$ satisfies $\phi(\alpha) \in \big(0, 1/{\tt h}\big) ~mod~ \bz$, then $\alpha$ is positive with respect to $B$. \eprop

\bsem Assumption \ref{basicassumptions} on the reductive group $G$ continues to be in force. Let $\mathfrak G$ be a subgroup scheme of $G$, and $H$ a maximal connected subgroup scheme of multiplicative type of $\mathfrak G$. For the existence of $H$ see \cite[Proposition 2.1]{deligne}.

Let $X(H)$ be the group of characters of $H$. The action of $H$ on ${\text {Lie}}(\mathfrak G)$ by conjugation gives an $X(H)$-gradation ${\text {Lie}}(\mathfrak G) = {\text {Lie}}(Z^{^0}_{_{\mathfrak G}}(H)) \oplus \bigoplus_{_{\alpha \neq 0}}{\text {Lie}}(\mathfrak G)^{^\alpha}$. Recall that if $M$ is a {\em connected} group scheme of multiplicative type over $k$, the group of characters $X(M)$ is finitely generated and does not have torsion prime to $p$. As $k$ is algebraically closed, this means that $M$ is a product of subgroups isomorphic to $\bg_{_m}$ or $\mu_{_{p^{^j}}}$ with $j > 0$.
\esem

\bcor\label{delignequestion} {\rm (cf. \cite[Lemma 2.12]{deligne})} Let $M$ be a connected subgroup scheme  of multiplicative type of  $G$ and suppose that $p~ >~{\tt h}_{_G}$.   The action of $M$ on ${\text {Lie}}(G)$ by conjugation gives an $X(M)$-gradation ${\text {Lie}}(G) =  \bigoplus_{_{\alpha \in X(M)}}{\text {Lie}}(G)^{^\alpha}$. If $\alpha \neq 0$, then for $\gamma \in {\text {Lie}}(G)^{^\alpha}$ one has $\gamma^{^{[p]}} = 0$. In particular, if $M = H$ and $\gamma \in {\text {Lie}}(\mathfrak G)^{^\alpha}$ we have $\gamma^{^{[p]}} = 0$. \ecor

\bpr If $\alpha$ is of infinite order, it suffices to apply \ref{2.10}. If $\alpha$ is of finite order, there is a subgroup $N$ of $M$ isomorphic to $\mu_{_{p^{^j}}}$ to which $\alpha$ restricts non-trivially and it suffices to prove \ref{delignequestion} for $N$ and the restriction of $\alpha$ to $N$. Let us choose an isomorphism of $X(N)$ with $X(\mu_{_{p^{^{j}}}})$ such that the image of $\alpha$ in $\bz/p^{^{j}}$ is of the form $p^{^{b}}$ with $0 \leq b < j$.

In a smooth algebraic group, the maximal connected subgroups of multiplicative type are (maximal) tori. Indeed, if $T$ is such a maximal subgroup, its centralizer $Z_{_{G}}(T)$ is smooth, being the fixed locus of a linearly reductive group acting on a smooth variety (see for example \cite[Theorem 2.8, Chapter II, \S5]{demgab}). Define $U:= Z_{_{G}}(T)^{^o}/T$. The group scheme $U$ is smooth, as a quotient of   $Z_{_{G}}(T)^{^o}$ which is smooth, and is unipotent by maximality of $T$. By \cite[Expos\'e XVII, Theorem 6.1.1]{sga3}, we have a splitting $Z_{_{G}}(T)^{^o} \simeq T \times U$, and hence $T$ is smooth.

In our case $N$ is hence contained in a maximal torus $T$ of $G$. Let 
\beqa
\phi:X(T) \to {\mathbb R}/{\bz}
\eeqa 
be the composite of 
$X(T) \to X(N)  =  \bz/p^{^{j}}$ and of the inclusion $x \mapsto x/{p^{^j}}$ of  $\bz/p^{^{j}}$ in ${\br}/{\bz}$.

The  $N$-weight space ${\text {Lie}}(G)^{^\alpha}$ is the sum of $T$-weight spaces ${\text {Lie}}(G)^{^\beta}$ for $\beta$ a root such that $\phi(\beta) = 1/p^{^{j-b}}$. As $\frac{1}{p^{^{j-b}}} \leq \frac{1}{p^{^{j}}} < \frac{1}{h}$ we may apply  \ref{yun}, to conclude that ${\text {Lie}}(G)^{^\alpha}$ is contained in the  in the Lie algebra of the unipotent radical of a Borel subgroup of $G$. In particular, it consists of nilpotent elements, proving \ref{delignequestion}. 
\epr

\bsem\label{completion}
The proof of \ref{strong}  now follows \cite[pages 594--599]{deligne} verbatim with a sole alteration; recall that in \cite{deligne} the group $G$ was the linear group  ${\text{GL}}(V)$ and the condition on the characteristic was $p > \dim(V)$. For an arbitrary connected reductive $G$, this condition now gets replaced by $p > {\tt h}_{_G}$, which makes \ref{delignequestion} applicable.\esem

\brem  If  $G = {\text{GL}}(V)$ one has ${\tt h}_{_G} = \dim(V)$. In the case $G = \prod {\text{GL}}(V_{_i})$ with for each $i$, $p > {\dim}(V_{_i})$, the case $p > {\tt h}_{_G}$ of   Theorem~\ref{strong} gives us  \cite[Theorem 1.7]{deligne}.  \erem

\begin{example}\rm (Brian Conrad) Here is an example in any characteristic $p > 0$, of a connected group of multiplicative type $M$ acting on a reductive group $G$, and of a non-trivial character $\alpha$ of $M$, such that the weight space ${\text {Lie}}(G)^\alpha$ contains elements which are not nilpotent. We take $G = {\text {SL}}_p$ (so that ${\tt h}_{_G} = p$) and  $M = \mu_{p^2}$, and for $\alpha$ the character $\zeta \mapsto \zeta^{^{[p]}}$.

We embed $M$ in the maximal torus of diagonal matrices of ${\text {SL}}_p$ by 
\[
\zeta \mapsto \mathrm{diag}(\zeta^{^0}, \zeta^{^{-p}}, \zeta^{^{-2p}}, \ldots, \zeta^{^{-(p-1)p}}).\] 
The restriction to $M$ of each simple root and of the lowest root is the character $\alpha: \zeta \mapsto \zeta^{^{[p]}}$, and ${\text {Lie}}(G)^\alpha$ is the sum of the corresponding root spaces. In the standard visualization of ${\text {SL}}_p$ this weight space inside ${\text {Lie}}(G) = \mathfrak {sl}_p$ is the span of the super-diagonal entries and the lower-left entry.

A sum of nonzero elements in those root lines contributing to ${\text {Lie}}(G)^\alpha$ is a  $p \times p$-matrix  $X \in\mathfrak {sl}_p$ which satisfies
\beqa
X(e_1) = t_p e_p, X(e_2) = t_1 e_1, \ldots, X(e_p) = t_{p-1} e_{p-1}.
\eeqa
Iterating $p$ times gives $X^p = \mathrm{diag}(t, \ldots, t)$, with $t := \prod t_j \ne 0$. Hence, $X \in {\text {Lie}}(G)^\alpha$ is not nilpotent. 
\end{example}

\begin{example}\label{deligneconrad} \rm A variant of the above example leads to an example of an infinitesimally saturated group scheme $\mathfrak G \subset SL(V)$ with $\dim(V) = p$ and such that $V$ is an irreducible representation and  $\mathfrak G_{_{red}}$ is a unipotent group. This in particular implies that $\mathfrak G_{_{red}}$ is not normal in $\mathfrak G$. 

Let $V$ be the affine algebra of $\mu_{_p}$, that is
\beqa
V=  \co(\mu_{_p}) := \frac{k[u]}{(u^{^{[p]}} - 1)}.
\eeqa
The vector space $V$ admits the basis $\{u^{^i} \mid i \in \bz/p\}$.

The multiplicative group $\co(\mu_{_p})^{^*}$ acts by multiplication on $V$. For $f \in \co(\mu_{_p})^{^*}$, $f^{^{[p]}}$ is constant. Define $Nf$ to be the constant value of $f^{^{[p]}}$. It is in fact the norm of $f$. The action of $\co(\mu_{_p})^{^*}$ on $V$ induces an action of the subgroup $N \subset \co(\mu_{_p})^{^*}$ for which $Nf = 1$.

On $V$ we have also the action of $\mu_{_p}$ by translations and this action normalizes the group $\co(\mu_{_p})^{^*}$ and its subgroup $N$.
Consider the group scheme:
\beqa
\mathfrak G := \mu_{_p} \ltimes N.
\eeqa
We make a few observations on $\mathfrak G $:
\begin{enumerate}[\rm (1)]
\item The group $\mathfrak G $ is infinitesimally saturated.
\item In $\mathfrak G $ consider the first factor $\mu_{_p}$ and  $\{u^{^i} \mid i \in \bz/p\} \simeq \bz/p$. These subgroups generate a subgroup $H$, which is a  central extension of $\mu_{_p}  \ltimes \bz/p$ by $\mu_{_p}$.
\item The representation $V$ is irreducible as an $H$-module. Indeed, a $\mu_{_p}$-submodule is generated by some of the $u^{^i}$, and if non-zero and $\bz/p$-stable, contains all of them. A fortiori, $V$ is an irreducible  $\mathfrak G$-module.
\item The reduced group $\mathfrak G _{_{red}}$ can be identified with the unipotent group $\{f \in \co(\mu_{_p})^{^*} \mid f(1) = 1\}$.
\item The subgroup $\mathfrak G _{_{red}}$ is not normal as  the point $1$ of $\mu_{_p}$ is not invariant by translations. 
\end{enumerate} 
\end{example}

\section{Semi-simplicity statements}

Let $G$ be a reductive group. Let $C$ be an algebraic group and $\rho:C \to G$ be a morphism.

\bdefe\label{gcr} (\cite[page 20]{mour}) One says that $\rho$ is {\em cr} if, whenever $\rho$ factors through a parabolic $P$ of $G$, it factors through a Levi subgroup of $P$. \edefe

When $G = {\text{GL}}(V)$, $\rho$ is {cr}  if and only if the representation $V$ of $C$ is completely reducible (or equivalently, semi-simple) and hence the terminology.

The property of $\rho$ being {cr} depends only on the subgroup scheme of $G$ which is the (schematic) image of $C$. It  in fact only depends on the image of $C$ in the adjoint group $G^{^{ad}}$. Indeed, the parabolic subgroups of $G$ are the inverse images of the parabolic subgroups of $G^{^{ad}}$, and similarly for the Levi subgroups. A subgroup scheme $\mathfrak G$ of $G$ will be called { cr} if its inclusion in $G$ is so.

For an irreducible root system $R$, let $\alpha_{_0}$ be the highest root and $\sum n_{_i} \alpha_{_i}$ its expression as a linear combination of the simple roots. The characteristic $p$ of $k$ is called {\em good} for $R$ if $p$ is larger than each $n_{_i}$. For a general root system $R$, $p$ is good if it is so for each irreducible component of $R$.

\bprop\label{likedeligne} Suppose $C$ is an extension 
\beqa
1 \to B \to C \to A \to 1
\eeqa
with $A^{^o}$ of multiplicative type and $A/A^{^o}$ a finite group of order prime to $p$, and suppose that $p$ is good. Let $\rho:C \to G$ be a morphism. If the restriction of $\rho$ to $B$ is {cr}, then $\rho$ is cr.\eprop

We don't know whether the proposition holds without the assumption that $p$ is good.
\bpr Let $P$ be a parabolic subgroup, $U$ its unipotent radical, and $\mathfrak u$ the Lie algebra of $U$. The parabolic  $P$ is said to be {\em restricted} if the nilpotence class of $U$ is less than $p$. If $P$ is a maximal parabolic corresponding to a simple root $\alpha$, the nilpotence class of $U$ is the coefficient of $\alpha$ in the highest root. It follows that $p$ is {\em good} if and only if all the maximal parabolic subgroups are restricted. By \cite[Proposition 5.3]{seitz} (credited by the author to Serre, see \ref{seitzetc}), if $P$ is restricted, one obtains by specialization from characteristic $0$, a $P$-equivariant isomorphism:
\beqa
\exp:(\mathfrak u, \circ) \stackrel{\sim}\to U
\eeqa
from $\mathfrak u$ endowed with the Campbell-Hausdorff group law, to $U$.

We will first show that whenever $\rho$ factors through a restricted parabolic subgroup $P$ as above, if the restriction of $\rho$ to $B$ factors through some Levi subgroup of $P$, the $\rho$ itself factors through some Levi subgroup of $P$.

The group $U(k)$ acts on the right on the set $\mathcal L(k)$ of Levi subgroups of $P$ by
\begin{center}
$u$ in $U(k)$ acts by $L \mapsto u^{^{-1}} L u$.
\end{center}
This action turns $\mathcal L(k)$ into a $U(k)$-torsor. This expresses the fact that two Levi subgroups are conjugate by a unique element of $U(k)$. The group $P(k)$ acts on $\mathcal L(k)$ and on $U(k)$ by conjugation. This turns $\mathcal L(k)$ into an equivariant $U(k)$-torsor.

We will need a scheme-theoretic version of the above. Fix a Levi subgroup $L_{_o}$. Let $\eL$ be the trivial $U$-torsor (i.e. $U$ with the right action of $U$ by right translations). We have the family of Levi subgroups $L_{_u}:= u^{^{-1}} L_{_o} u$ parametrized by $\eL = U$. We let $P$ act on $U$ by conjugation, and on $\eL$ as follows: $p = v \ell$ in $P = U L_{_o}$ acts on $\eL = U$ by  $u \mapsto v^{^{-1}}.pup^{^{-1}}$. This turns $\eL$ into an equivariant $U$-torsor. When we pass to $k$-points, and attach to $u$ in $\eL$ the Levi subgroup $u L_{_o} u^{^{-1}}$, we recover the previously described situation.

The morphism $\rho:C \to P$ turns $\eL$  into an equivariant $U$-torsor. A point $x$ of $\eL$ corresponding to a Levi subgroup $L_{_x}$ is fixed by $C$ (scheme-theoretically) if and only if $\rho$ factors through $L_{_x}$. This expresses the fact that a Levi subgroup is its own normalizer in $P$.

We want to prove that if $B$ has a fixed point in $\eL$, so does $C$. Let $U^{^B}$ be the subgroup of $U$ fixed by $B$, for the conjugation action. If $B$ has a fixed point in $\eL$, the fixed locus $\eL^{^B}$ is a $U^{^B}$-torsor. As $B$ is a normal subgroup of $C$, $C$ acts on $\eL^{^B}$ and $U^{^B}$, and the action factors through $A$.

The isomorphism ${\text{exp}}:(\mathfrak u, \circ) \to U$ is compatible with the action of $B$ by conjugation. Hence it induces  an isomorphism from $(\mathfrak u^{^B}, \circ) \to U^{^B}$. Let $Z^{^i}(\mathfrak u^{^B})$ be the central series of $\mathfrak u^{^B}$, and define $Z^{^i}(U^{^B}):= {\text{exp}}(Z^{^i}(\mathfrak u^{^B}))$. The isomorphism ${\text{exp}}$ induces an isomorphism between the vector group $Gr^{^i}_{_Z}(\mathfrak u^{^B})$ and $Z^{^i}(U^{^B})/Z^{^{i+1}}(U^{^B})$, compatible with the action of $A$. On $Gr^{^i}_{_Z}(\mathfrak u^{^B})$ , this action is linear.

The assumption on $A$ amounts to saying that $A$ is linearly reductive, that is, all its representations are semi-simple. Equivalently, if $k$ is the trivial representation, any extension
\beqa
0 \to V \stackrel{a}\to E \stackrel{b}\to k \to 0
\eeqa
splits. Passing from $E$ to $b^{^{-1}}(1)$, such extensions correspond to $A$-equivariant $V$-torsors, and the extension splits if and only if $A$ has a scheme-theoretic fixed point on the corresponding torsor.

Define $U_{_i}^{^B}$ to be $U^{^B}/Z^{^i}(U^{^B})$, and $\eL_{_i}^{^B}$ to be the $U_{_i}^{^B}$-torsor obtained from $\eL$ by pushing by $U^{^B} \to U_{_i}^{^B}$. We prove by induction on $i$ that $A$ has a fixed point on $\eL_{_i}^{^B}$.

As $U_{_1}^{^B}$ is trivial, the case $i =1$ is trivial. If $x$ is a fixed point of $A$ in $\eL_{_i}^{^B}$, the inverse image of $x$ in $\eL_{_{i+1}}^{^B}$ is an equivariant $A$-torsor on 
$Gr^{^i}_{_Z}(\mathfrak u^{^B}) \sim Z^{^i}(U^{^B})/Z^{^{i+1}}(U^{^B})$. By linear reductivity, $A$ has a fixed point on the inverse image. As the central descending series of $\mathfrak u$, and hence of $\mathfrak u^{^B}$ terminates, this proves \ref{likedeligne} for restricted parabolic subgroups. 

We now prove \ref{likedeligne} by induction on (the dimension) of $G$. Suppose that $\rho$ factors through a proper parabolic subgroup $P$. As $p$ is {\em good}, there exists a restricted proper parabolic $Q$ containing $P$, and $P$ is the inverse image by the projection $Q \to Q/R_{_u}(Q)$ of a parabolic $P'$ of $Q/R_{_u}(Q)$. Let $L$ be a Levi subgroup of $Q$ through which $\rho$ factors and let $P'_{_L}$ be the parabolic subgroup of $L$ obtained as the inverse image of $P'$ by the isomorphism $L \stackrel{\sim}\to Q/R_{_u}(Q)$. Levi subgroups of $P'_{_L}$ are Levi subgroups of $P$, and it remains to apply the induction hypothesis to $L$, for which $p$ is good too.
\epr

\brem For several results related to \ref{likedeligne} but in the setting of reduced subgroups, see \cite[Theorem 3.10]{bmr1} and \cite[Theorem 1.1 and Corollary 3.7]{bmr2}.
\erem

 Fix in the reductive group $G$, a maximal torus $T$, and a system of simple roots corresponding to a Borel subgroup $B$ containing $T$. Let $U$ be the unipotent radical of $B$. 

\bdefe\label{height}(cf. \cite{dynk}, \cite{serre2}, \cite{imp}) 
 The
Dynkin {\em height}  $ht_{_G}(V)$ of a representation $V$ of $G$  is the largest among  
$\{\sum_{_{\alpha > 0}} \langle \lambda, \alpha^{^\vee}\rangle\}$, for $\lambda$ a weight for the action of $T$ on
$V$. \edefe

 This notion and this terminology go back to Dynkin \cite[pages 331--332]{dynk} where it is called ``height" as in \cite{imp}, while in Serre (\cite{serre2} and \cite{mour}), it is simply ``$n(V)$".  If $V$ is an irreducible representation, with dominant weight $\lambda^{^+}$ and smallest weight $\lambda^{^{-}}$, it is the sum of the coefficients of $\lambda^{^+} - \lambda^{^{-}}$, expressed as linear combination of the simple roots. 

It follows that the product in ${\text{End}}(V)$ of the action of $ht_{_G}(V) + 1$  elements of ${\text{Lie}}(U)$ vanishes and for  ${\text{\cursive n}}$ nilpotent  in ${\text {Lie}}(G)$, one has  ${\text{\cursive n}} ^{^{ht(V) + 1}} = 0$ in ${\text{End}}(V)$. 
 
\bsem\label{lowheight} The representation $\rho:G \to {\text{GL}}(V)$ is said to be of {\em low
    height} if $p>ht_{_G}(V)$.  By \cite[Theorem 6, page 25]{mour}, representations of low height are semi-simple. One can show that if $G$ admits a representation $V$ of low height  which is {\em almost faithful}, meaning that its kernel  is of multiplicative type, then $G$ satisfies the assumption \ref{basicassumptions}. That $p \geq {\tt h}_{_G}$ results from the more precise statement that $ht_{_G}(V) \geq {\tt h}_{_G} - 1$ (\cite[(5.2.4), page 213]{serre3}). For the property that $\widetilde{G}/G'$ is \'etale, one uses the fact that the non-trivial  irreducible representation of ${\text{PGL}}(p)$ of the smallest height is the adjoint representation which is of height $2p-2$.
    
 We now assume that $V$ is of low height, and the assumption \ref{basicassumptions} holds for $G$. It follows that any nilpotent    ${\text{\cursive n}}$ in ${\text{Lie}}(G)$ satisfies ${\text{\cursive n}} ^{^{[p]}} = 0$, and further, the exponential map \eqref{springexp} is defined. The image $d\rho({\text{\cursive n}})$ of ${\text{\cursive n}}$ in ${\text{Lie}}(\text{GL}(V)) = \text{End}(V)$ also has a vanishing $p^{^{th}}$-power, hence $\text{exp}(d\rho({\text{\cursive n}}).t)$ is defined. By \cite[Theorem 5, page 24]{mour}, one has the following compatibility statement.
 \esem

\bsem\label{compatibilityofspring} \underline{Compatibility}:  If ${\text{\cursive n}}$  in ${\text{Lie}}(G)$  is nilpotent 
\beqa\label{compatibility2}
\rho({\text {exp}}(t {\text{\cursive n}})) = {\text {exp}}(t d\rho({\text{\cursive n}})).
\eeqa
As a consequence, if $u^{^{[p]}} = 1$ in $G$, one has
\beqa\label{compatibility}
\rho(u^{^t}) = \rho(u)^{^t}.
\eeqa
\esem

The following theorem is a schematic analogue of \cite{bt}, for $p$ large enough.
\bth\label{schematicbt} Suppose that the reductive group $G$ admits a low height almost faithful representation $\rho:G \to {\text{GL}}(V)$, and that $p > {\tt h}_{_G}$. Then, for any non-trivial unipotent subgroup $U$ of $G$, there exists a proper parabolic subgroup $P$ of $G$ containing the normalizer $N_{_G}(U)$ of $U$, and whose unipotent radical contains $U$. \eeth

The condition $p > {\tt h}_{_G}$ implies that $G$ satisfies the assumption \ref{basicassumptions}. If $G$ is simple simply connected, it implies the existence of an almost faithful low height representation except for the $G$ of type $F_{_4}$, $E_{_6}$, $E_{_7}$ or $E_{_8}$, in which case the lowest height of a non-trivial representation and the Coxeter number are respectively $16 > 12$, $16 > 12$, $27 > 18$ and $58 > 30$. For these groups, we do not know whether the conclusion of the theorem is valid assuming only $p > {\tt h}_{_G}$. 

\bpr Let $V^{^U}$ be the invariants of $U$ acting on $V$. It is not zero, because $U$ is unipotent. It is not $V$, because the representation $V$ is almost faithful, hence faithful on $U$. It does not have a $U$-stable supplement $V'$ in $V$, because $U$ would have invariants in $V'$.

Let $H$ be the subgroup scheme of $G$ which stabilizes $V^{^U}$. It contains the normalizer $N_{_G}(U)$ of $U$. It is a doubly saturated subgroup scheme of $G$. Indeed, if $h$ in $H(k)$ is of order $p$, by \eqref{compatibility}, $\rho(h^{^t}) = \rho(h)^{^t} = \sum_{_{i < p}} \binom{t}{i}(\rho(h)-1)^{^i}$, which stabilizes $V^{^U}$, and similarly if ${\text{\cursive n}}$ in $\text{Lie}(H)$ is nilpotent, the $\text{exp}(t {\text{\cursive n}})$ are in $H$.

As $p > {\tt h}_{_G}$, theorem \ref{strong} ensures that $H_{_{red}}^{^0}$ is a normal subgroup scheme of $H$ and that the quotient $H/H_{_{red}}^{^0}$ is an extension of a finite group of order prime to $p$ by a group of multiplicative type. It follows that $U \subset H_{_{red}}^{^0}$.

\blem\label{4.6} $H_{_{red}}^{^0}$ is not reductive. \elem

\bpr If it were, $V$ would be a representation of low height of $H_{_{red}}^{^0}$ (\cite[Corollary 1, page 25]{mour}), hence a semi-simple representation of $H_{_{red}}^{^0}$, and $V^{^U}$ would have in $V$ a $H_{_{red}}^{^0}$-stable supplement. As $V^{^U}$ does not admit a supplement stable under $U \subset H_{_{red}}^{^0}$, this is absurd.\epr

\noindent {\underline{Proof of \ref{schematicbt} continued}}: If $S$ is a doubly saturated subgroup scheme of $G$, we will call $R_{_u}(S_{_{red}}^{^0})$ the {\em unipotent radical} of $S$ and denote it simply by $R_{_u}(S)$. By \ref{strong}, it is a normal subgroup of $S$ and $S/R_{_u}(S_{_{red}}^{^0})$ does not contain any normal unipotent subgroup. This justifies the terminology.

Define $U_{_1} := R_{_u}(H)$. By \ref{4.6}, it is a non-trivial unipotent subgroup of $G$, and we can iterate the construction. We define for $i \geq 1$
\begin{center}
$H_{_i}$:=  stabilizer of $V^{^{U_{_{i}}}} \subset V$

$U_{_{i+1}}:= R_{_u}(H_{_i})$.
\end{center}
One has $U \subset H_{_{red}}^{^0}$ and 
\beqa
U \subset N_{_G}(U) \subset H \subset N_{_G}(U_{_1}) \subset H_{_1} \subset N_{_G}(U_{_2}) \subset H_{_2} \ldots
\eeqa
The $H_{_{i,red}}^{^o}$ form an increasing sequence of smooth connected subgroups of $G$. It stablizes, hence so do the sequences of the $U_{_i}$ and of the $H_{_i}$. If $H_{_i} = H_{_{i+1}}$, one has $H_{_i} = N_{_G}(U_{_{i+1}}) = H_{_{i+1}}$ and 
\begin{center}
$U_{_{i+1}} = R_{_u}(H_{_i}) = R_{_u}(H_{_{i+1}}) = R_{_u}(N(U_{_{i+1}})_{_{red}}^{^o})$.
\end{center}
By \cite[Proposition 2.3, page 99]{bt} (or for example \cite[Section 30.3, Proposition on page 186]{humphreys}), this implies that $N(U_{_{i+1}})_{_{red}} = (H_{_{i+1}})_{_{ red}}$ is a proper parabolic subgroup of $G$. Call it $Q$. A parabolic subgroup of $G$ is its own normalizer scheme (cf. \cite[XII, 7.9]{sga3}, \cite[page 469]{cgp}). As $(H_{_{i+1}})_{_{ red}} = Q$ is normal in $H_{_{i+1}}$, it follows that $H_{_{i+1}} = Q$ and that $N_{_G}(U) \subset Q$.

A Levi subgroup $L$ of $Q$ is a reductive subgroup of $G$, hence satisfies the assumptions of \ref{schematicbt}. It is isomorphic to $Q/R_{_u}(Q)$. If $U$ is not contained in $R_{_u}(Q)$, we can repeat the argument for the image $\overline{U}$ of $U$ in $Q/R_{_u}(Q)$, which is isomorphic to $L$. One obtains a proper parabolic subgroup of $Q/R_{_u}(Q)$ which contains the normalizer of $\overline{U}$. Its inverse image in $Q$ is a parabolic subgroup, properly contained in $Q$ and containing the normalizer of $U$. Iterating, one eventually finds a parabolic $P$ containing $N_{_G}(U)$ and such that $U \subset R_{_u}(P)$.
\epr

\brem\label{4.7} Let $\widetilde{G}$ be the simply connected central extension of the derived group $G'$, and let $\widetilde{G}_{_i}$ be its simple factors: $\widetilde{G} = \prod \widetilde{G}_{_i}$. A representation of low height $V$ of $G$ is almost faithful if and only if its restriction to each $\widetilde{G}_{_i}$ is not trivial. It suffices to check this for each $\widetilde{G}_{_i}$ separately. Thus we may assume $G$ is simply connected. The existence of a non-trivial $V$ of low height implies that $p > 2$ for $G$ of type $B_{_n}$, $C_{_n} (n \geq 2)$ or $F_{_4}$ and $p > 3$ for $G$ of type $G_{_2}$. Let $\overline{G}$ be the image of $G$ in $\text{GL}(V)$. If $V$ is non-trivial, $u:G \to \overline{G}$ is an isogeny. We want to show that it is a central isogeny. If it is not, the structure of isogenies (\cite[XXII, 4.2.13]{sga3}) shows that $\ker(u)$ contains the kernel of the Frobenius. The weights of $V$ are the $p^{^{th}}$-powers and $ht_{_G}(V)$ is a multiple of $p$, contradicting the low height assumption. \erem 

\bsem If $G$ is a reductive group for which assumption \ref{basicassumptions} holds, and $\mathfrak G$ is a subgroup scheme of $G$, the {\em double saturation}  $\mathfrak G^{^*}$ of $\mathfrak G$ is the smallest doubly saturated subgroup scheme of $G$ containing $\mathfrak G$. It is the intersection of the doubly saturated subgroup schemes containing $\mathfrak G$, and is obtained from $\mathfrak G$ by iterating the construction of taking the group generated by $\mathfrak G$, the additive groups $\text{exp}(t {\text{\cursive n}})$ for ${\text{\cursive n}}$ nilpotent in $\text{Lie}(\mathfrak G)$ and $u^{^t}$ for $u$ of order $p$ in $\mathfrak G(k)$.\esem

\bcor\label{4.8} Let $V$ be a low height almost faithful representation of a reductive group $G$. Assume that $p > {\tt h}_{_G}$. Let $\mathfrak G$ be a subgroup scheme of $G$, and let $\mathfrak G^{^*}$ be its double saturation. Then the following conditions are equivalent:
\begin{enumerate}[\rm (i)]
\item $V$ is a semi-simple representation of $\mathfrak G$
\item $V$ is a semi-simple representation of $\mathfrak G^{^*}$
\item $\mathfrak G$ is cr in $G$
\item $\mathfrak G^{^*}$ is cr in $G$
\item the unipotent radical of $(\mathfrak G^{^*})_{_{red}}^{^0}$ is trivial.
\end{enumerate}
\ecor

\bpr (i) $\iff$ (ii): If $W$ is a subspace of $V$, the stabilizer in $G$ of $W$ is doubly saturated, as we saw in the beginning of the proof of \ref{schematicbt}. If $\mathfrak G$ stabilizes $W$, it follows that $\mathfrak G^{^*}$ also stabilizes $W$: the lattice of sub-representations of $V$ is the same for $\mathfrak G$ and $\mathfrak G^{^*}$, hence the claim.

(iii) $\iff$ (iv): Similarly, the parabolic subgroups of $G$ and their Levi subgroups are doubly saturated, hence contain $\mathfrak G$ if and only if they contain $\mathfrak G^{^*}$.

not (v) $\implies$ not (ii): Let $U$ be the unipotent radical of $(\mathfrak G^{^*})_{_{red}}^{^0}$. If it is non-trivial, $V^{^U} \neq V$, because $V$ is faithful on $U$ and does not have a $U$-stable supplement. As $U$ is normal in $\mathfrak G^{^*}$ (\ref{strong}), $V^{^U}$ is a sub-representation for the action of $\mathfrak G^{^*}$ on $V$. This contradicts the semi-simplicity of $V$.

(v) $\implies$ (ii): The representation $V$ of the reductive group $(\mathfrak G^{^*})_{_{red}}^{^0}$ is of low height, hence semi-simple. By \ref{strong}, $(\mathfrak G^{^*})_{_{red}}^{^0}$ is a normal subgroup of $\mathfrak G^{^*}$ and the quotient $A$ is linearly reductive. If $W$ is a sub-$\mathfrak G^{^*}$-representation of $V$, $A$ acts on the affine space of $(\mathfrak G^{^*})_{_{red}}^{^0}$-invariant retractions $V \to W$. It has a fixed point, whose kernel is a supplement to $W$.

not (v) $\implies$ not (iv): Let $U$ be the unipotent radical of $(\mathfrak G^{^*})_{_{red}}^{^0}$.  If it is not trivial, there exists a parabolic $P$ containing its normalizer, hence $\mathfrak G^{^*}$, and the unipotent radical of $P$ contains $U$ (\ref{schematicbt}). Thus, no Levi subgroup of $P$ can contain $\mathfrak G^{^*}$.

(v) $\implies$ (iv): By \cite[Theorem 7, page 26]{mour}, $(\mathfrak G^{^*})_{_{red}}^{^0}$ is cr in $G$, and one applies \ref{likedeligne}.
\epr

\bcor\label{delignestr3}  Let $\rho:G \to {\text{GL}}(V)$ be an almost faithful  low height representation, and let $v$ in $V$ be an element such that the $G$-orbit of $v$ in $V$ is {\em closed}. Then there exists a connected multiplicative central subgroup scheme $M \subset G^{^0}_{_v}$ and a surjective homomorphism $M \times G^{^0}_{_{v,red}} \to G^{^0}_{_v}$. \ecor

\bpr The orbit being closed in $V$ and hence affine, the reduced stabilizer  $G_{_{v,red}} $ is reductive (\cite{borel}). Since $p \geq {\tt h}_{_G}$ and  since  stabilizers are doubly saturated and  \ref{strong}(2) holds for $\mathfrak G = G_{_v}$,  we get the required result. 
\epr

We now observe that the results of \cite[Section 6]{deligne} can be obtained as a consequence of  \ref{4.8}.  Note that by the remarks in \cite[page 607]{deligne} it suffices to prove the semisimplicity results in the case when $k$ is algebraically closed. 

\bth\label{section6} $\mathfrak G$ be  an algebraic group. Let $(V_{_i})_{_{i \in I}}$ be a finite family of semi-simple $\mathfrak G$-modules  and let $m_{_i}$ be integers $\geq 0$. If 
\beqa\label{extlowht}
\sum m_{_{i}} (\dim~V_{_{i}} - m_{_{i}}) < p
\eeqa
the $\mathfrak G$-module $\bigotimes_{_{j}} \bigwedge^{^{m_{_j}}} V_{_{j}}$ is semi-simple. \eeth

\bpr Let $G = \prod_{_{j}} {\text{GL}}(V_{_{j}})$ and $V = \bigotimes_{_{j}} \bigwedge^{^{m_{_j}}} V_{_{j}}$. Then, \eqref{extlowht} is simply the inequality $p > ht_{_{G}}(V)$. Replacing $\mathfrak G$ by its image in $G$ we may and shall assume that $\mathfrak G$ is a subgroup scheme of $G$. Since $V_{_i}$ are semi-simple $\mathfrak G$-modules, it follows that $\mathfrak G$ is $G$-cr, for $G = \prod_{_{j}} {\text{GL}}(V_{_{j}})$. By \cite[\S 6.2]{deligne}, we may also assume that $p > {\dim}(V_{_j}) = {\tt h}_{_{{\text{GL}}(V_{_{j}})}} \forall j$.

Hence by working with the image of $G$ (and $\mathfrak G$) in ${\text{GL}}(V)$ and applying \ref{4.8}, we conclude that $V$ is semi-simple as a $\mathfrak G$-module.
\epr

\brem If $(V_{_i}, q_{_i})$ is a  non-degenerate quadratic space with $\dim~V_{_i} = 2 d_{_i}$ on which $\mathfrak G$ acts by similitudes, then by passing to a subgroup of index at most $2$ and mapping to the group of similitudes rather than $\text{GL}(V_{_i})$, one can replace  the term  $m_{_{i}} (\dim~V_{_{i}} - m_{_{i}})$ by $m_{_{i}} (\dim~V_{_{i}} - m_{_{i}} - 1)$, when $m_{_i} \leq 2d_{_i}$. \erem

{\em Complete reducibility in the classical case}.
By \ref{gcr}, a subgroup scheme $\mathfrak G \subset G$ is called cr if for every parabolic subgroup $P \subset G$ containing $\mathfrak G$, there exists an opposite parabolic subgroup  $P'$ such that $\mathfrak G \subset P \cap P'$. Suppose that $char(k) \neq 2$ and let $G$ be $SO(V)$ (or $Sp(V)$ in any characteristic), relative to a non-degenerate symmetric or alternating bilinear form $B$ on $V$. In this situation, the notion of cr can be interpreted  as follows: $\mathfrak G$ is cr in $G$ if and only if for every $\mathfrak G$-submodule $W \subset V$ which is totally isotropic, there exists a totally isotropic $\mathfrak G$-submodule  $W'$ of the same dimension, such that the restriction of $B$ to $W + W'$ is non-degenerate (cf. \cite[Example 3.3.3, page 206]{serre3}).

\blem\label{serreselementaryargument}  Let the subgroup scheme $\mathfrak G$ of  $G$ be cr.  Then the $\mathfrak G$-module $V$ is semi-simple and conversely. \elem

\bpr  Let $W \subset V$ be a $\mathfrak G$-submodule. Then we need to produce a $\mathfrak G$-complement.

Consider $W_{_1}:= W \cap W^{^\perp}$. If $W_{_1} = (0)$, then $W \oplus W^{^\perp} = V$ and we are done. So let $W_{_1} \neq (0)$.  Then $W_{_1}$ is a $\mathfrak G$-submodule which is totally isotropic and hence by the cr property, we have a totally isotropic $\mathfrak G$-submodule $W_{_1}'$ of the same dimension as $W_{_1}$, such that the form $B$ is non-degenerate on $W_{_1} + W_{_1}'$. In particular, $W_{_1} \cap W_{_1}' = (0)$. Since $W_{_1} \subset W^{^\perp}$, we see that $W \subset W_{_1}^{^\perp}$. 

Let $w \in W \cap W_{_1}' \subset W_{_1}^{^\perp} \cap W_{_1}'$ and suppose $w \neq 0$. Since $w \in W_{_1}'$, there exists $w' \in W_{_1}$ such that $B(w, w') \neq 0$. On the other hand, since $w \in W_{_1}^{^\perp}$, $B(w,v) = 0$ for all $v \in W_{_1}$ and in particular $B(w,w') = 0$ which contradicts the assumption that $w \neq 0$. Hence it follows that $W \cap W_{_1}' = (0)$. Thus, $W \subsetneq X = W \oplus W_{_1}' \subset V$ is a $\mathfrak G$-submodule.  We proceed similarly and get $X_{_1} = X \cap X^{^\perp}$ such that $X \oplus X_{_1}' \subset V$. If $X_{_1} = (0)$, then $V = X \oplus X^{^\perp}$, so get a $\mathfrak G$-decomposition of $V$ as $W \oplus W_{_1}' \oplus \ldots$, that is a $\mathfrak G$-complement of $W$ in $V$.  

Conversely, let $V$ be semi-simple as a $\mathfrak G$-module. Let $W \subset V$ be a totally isotropic $\mathfrak G$-submodule of $\dim(W) = d$. Note that $d < \dim(V)/2$. We therefore have a $\mathfrak G$-submodule $Z \subset V$ such $V = W \oplus Z$. The non-degenerate form $B$ gives an $\mathfrak G$-equivariant isomorphism $\phi:V \to V^{^*} = W^{^*} \oplus Z^{^*}$ and since $W$ is totally isotropic $\phi(W) \cap W^{^*} = (0)$. Hence $\phi(W) \subset Z$. 

Again, since $W$ is totally isotropic, the restriction of $B$ to $Z$ is non-degenerate and hence we get an isomorphism $\psi:Z^{^*} \to Z$. Define $W':= \psi \circ \phi(W)$.

Then it is easily seen that $W'$ is of dimension $d$ and also totally isotropic $\mathfrak G$-invariant submodule of $V$. 
Finally, $B$ is non-degenerate on $W \oplus W'$. 
Hence $\mathfrak G$ is cr in $G$. 
\epr

\section{\'Etale slices in positive characteristics}

Let $G$  act on an affine variety $X$ and let  $x$ be a point of $X$ whose orbit $G.x$ is closed in $X$. We prove a schematic analogue of Luna's \'etale slice theorem ( \ref{luna}, \ref{5.6}) for suitable bounds on the characteristic $p$. In order to make this precise, one needs to extend the notion of a slice in  \cite[Definition 7.1]{bardsley}  to a scheme theoretic setting (see \ref{5.6})  and in the process one has to work with schematic  stabilizers of closed orbits. The key issue then becomes the semi-simplicity of the tangent space to $x$ in $X$ as a representation of the stabilizer of $G$ at $x$. Here the structure theorem \ref{strong} become crucial and the bounds on $p$ are forced as a consequence.  In \cite{bardsley} the scheme-theoretic aspects are eschewed forcing them to make the assumption that  the orbit $G.x$ is ``separable'', i.e. the stabilizer $G_{_x}$ is reduced. The bounds on the characteristic therefore do not show up in \cite{bardsley}.

We begin this section with the following (linear) analogue of the Luna \'etale slice theorem in positive characteristics.
 
\bth\label{luna}
 Let $V$ be a $G$-module with low height i.e. such that $p~> ~ht_{_G}(V)$. Let $v$ in $V$ be an element such that the orbit $G.v$  is a closed orbit in $V$. Then there exists a $G_{_v}$-invariant linear subspace $S$ of $V$ giving rise to  a commutative diagram:
\beqa\label{cartesian}
\xymatrix{
G \times^{^{G_{_v}}} S  \ar[r]^-{\phi} \ar[d]_{f} &
V \ar[d]^{q} \\
(G \times^{^{G_{_v}}} S)\parallelslant G \ar[r]^-{\ell} &  V\parallelslant G
}
\eeqa 
and $G$-equivariant open subsets $U \subset (G \times^{^{G_{_v}}} S)$ containing the closed orbit $G.v$  and an open subset $U'$ of  $V$ containing $v$, for which \eqref{cartesian} induces a cartesian diagram 
\beqa\label{cartesian1}
\xymatrix{
U  \ar[r]^{\phi} \ar[d]_{f} &
U' \ar[d]^{q} \\
U\parallelslant G  \ar[r]^{\ell|_{_{U}}} &  U'\parallelslant G
}
\eeqa
such that the morphism $\ell|_{_{U}}$ is  \'etale.
\eeth

\brem The above theorem was stated in the note \cite{paramvikram}  whose proof contained serious gaps (as was pointed to the authors by Serre in a private correspondence). \erem

We note that $G/G_{_v}$ is constructed in \cite[III, Proposition 3.5.2]{demgab}. It represents the quotient in the category of fppf sheaves.  Furthermore if $\pi_{_v}:G  \to V,~ g \mapsto g.v$, the image $\text{im}(\pi_{_v})$,   as a locally closed sub-scheme of $V$ with its {\em reduced scheme} structure, can be identified with the scheme $G/G_{_v}$. We call this locally closed sub-scheme, the {\em orbit} $G.v$ and have the identification $G/G_{_v} \simeq G.v$ (see also \cite[Proposition and Definition  1.6, III, \S3, page 325]{demgab}).

\bprop\label{genpropo} Let $V$ be an arbitrary $G$-module. Let $v \in V$ and suppose that there exists a $G_{_v}$-submodule $S$ of the tangent space $T_{_v}(V)$, such that $T_{_v}(V)$ splits as $T_{_v}(G.v) \oplus S$. Let $G$ act on $G \times S$ by $h.(g,s):= (h.g,s)$. Then the $G$-morphism $\Phi:G \times S \to V$ given by $(g,s) \mapsto g.v + g.s$  descends to a $G$-morphism 
\beqa\label{phi}
\phi:G \times^{^{G_{_v}}} S \to V
\eeqa
which is {\em \'etale} at  $(e, 0)$, $e$ being the identity of $G$.
\eprop

\bpr To ensure that $\Phi$ descends to a morphism $\phi$,  we need to check for every commutative $k$-algebra $A$, that  $\Phi$ is constant on all $G_{_v}(A)$-orbits. For simplicity of exposition we will suppress the $A$. 

Let $\alpha$ in $G_{_v}$ act on $G \times S$ by $\alpha.(g,s) = (g.\alpha,\alpha^{^{-1}}.s).$ Observe that $\Phi\bigl(\alpha.(g,s)\bigr) = \Phi(g.\alpha,\alpha^{^{-1}}.s) = g.\alpha.v + g.\alpha.\alpha^{^{-1}}.s = g.v + g.s$ (since $\alpha$ fixes $v$). Therefore it is constant on the $G_{_v}$-orbits. Since the action of $G_{_v}$ on $G \times S$ is scheme-theoretically free, $\Phi$ descends to a morphism $\phi:G \times ^{^{G_{_v}}} S \to V$. Clearly the actions of $G$ and $G_{_v}$ on $G \times S$ commute and hence the descended morphism is also a $G$-morphism.

Observe that the quotient morphism $G \to G/G_{_v}$ is a torsor for the group scheme $G_{_v}$, locally trivial under the fppf topology.  Since the action of $G_{_v}$ on $S$ is linear, we see that the associated fibre space $\psi: G \times^{^{G_{_v}}} S \to G/G_{_v}$
is a locally free sheaf of rank $= \dim(S)$. In particular, $G \times^{^{G_{_v}}} S$ is a smooth $k$-scheme of finite type.  Observe further that under the morphism $\phi$, the zero section of the vector bundle $\psi: G \times^{^{G_{_v}}} S \to G/G_{_v}$ canonically maps onto the orbit $G.v \subset V$, while the fibre $~\psi^{^{-1}}(e.G_{_v})$ of the identity coset $e.G_{_v} \in G/G_{_v}$ maps isomorphically to the affine subspace $S + v \subset V$. Since $T_{_v}(V) = T_{_v}(G.v) \oplus S$  by assumption, it follows that the differential $d\phi_{_z}$, at $z = (e, 0)$, is an isomorphism.

We now apply \cite[Lemme 2.9]{deligne}, to conclude that the morphism $\phi~:~G \times^{^{G_{_v}}}~S~\to~V$ is \'etale at $z = (e, 0)$.
\epr

\bprop\label{lowhtimpliescr} Let $V$ be a $G$-module such that $p~> ~ht_{_G}(V)$, and let $v$ in $V$ be an
element such that the $G$-orbit of $v$ in $V$ is closed.  Then there exists a $G_{_v}$-submodule $S \subset V$ such that $V = T_{_v}(G.v) \oplus S$ as a $G_{_v}$-module. In particular, the consequences of \ref{genpropo} holds.  \eprop

\bpr One knows that $G_{_v}^{^0}$ is the kernel of the natural map $G_{_v} \to \pi_{_0}(G_{_v})$, to the group of connected components, and hence is a normal subgroup of finite index. Further, we note that $\mid~G_{_v}/G^{^0}_{_v} \mid\ =\ \mid ~G_{_{v,red}}/G^{^0}_{_{v,red}} \mid$.  Note also that since $G_{_{v,red}}$ is a saturated subgroup of $G$, by \cite[page 23, Property 3]{mour} the index $\mid ~G_{_{v,red}}/G^{^0}_{_{v,red}} \mid$ is prime to $p$.

Since  $G.v$ is a closed orbit by \ref{delignestr3} (and  \ref{strong}) we have an exact sequence:
\beqa
1 \to G^{^0}_{_{v,red}} \to G^{^0}_{_v} \to \tau \to 1
\eeqa
where $\tau$ is a multiplicative group scheme. 

Now $V$ is semi-simple as a $G^{^0}_{_{v,red}}$-module, $G^{^0}_{_{v,red}}$ being a saturated subgroup of $G$ and also as a $\tau$-module since it is  multiplicative (and hence linearly reductive). Thus by \cite[Lemme 4.2]{deligne} 
$V$ is semi-simple as a $G^{^0}_{_v}$-module and therefore as a $G_{_v}$-module as well.

In particular, we have a $G_{_v}$-supplement $S$ for the $G_{_v}$-invariant subspace $T_{_v}(G.v) \subset  T_{_v}(V) = V$, i.e. we have a $G_{_v}$-decomposition $S \oplus T_{_v}(G.v)$  for $V$.
\epr

We recall the ``Fundamental lemma of Luna" which holds in positive characteristics as well and which is essential to complete the proof of Theorem \ref{luna}.

\blem\label{fundlem} {\rm (\cite[page 152]{git})} Let $\text{\cursive s}:X \to Y$ be a $G$-morphism of affine $G$-schemes. Let $F \subset X$ be a closed orbit such that:
\begin{enumerate}[\rm (1)]
\item $\text{\cursive s}$ is \'etale at some point of $F$
\item $\text{\cursive s}(F)$ is closed in $Y$
\item $\text{\cursive s}$ is injective on $F$
\item $X$ is normal along $F.$
\end{enumerate}
Then there are affine $G$-invariant open subsets $U \subset X$ and $U' \subset Y$ with $F \subset U$ such that $\text{\cursive s}\parallelslant G:U \parallelslant G \to U'\parallelslant G$ is \'etale and $(\text{\cursive s}, p_{_U}): U \to U' \times_{_{U' \parallelslant G}} U \parallelslant G$ is an isomorphism (where $p_{_U}:U \to U \parallelslant G$ is the quotient morphism). \elem
 
\begin{proof}[Proof of  \ref{luna}] As a first step, we need to show that the categorical quotient $G~\times^{^{G_{_v}}}~S~\parallelslant G$ exists as a scheme.  As we have seen in  \ref{genpropo}, action of $G$ on $G \times S$ (via $g.(a,s) = (g.a,s)$ and the (twisted) action of $G_{_v}$ on $G \times S$ commute and hence
\beqa
k[G \times^{^{G_{_v}}} S]^{^G} = (k[G \times S]^{^{G_{_v}}})^{^G} =  (k[G \times S]^{^G})^{^{G_{_v}}} = k[S]^{^{G_{_v}}}.
\eeqa
Thus it is enough to show that $k[S]^{^{G_{_v}}}$ is finitely generated and we would have
\beqa
G \times^{^{G_{_v}}} S \parallelslant G \simeq S \parallelslant G_{_v}.
\eeqa

As we have observed in the proof of  \ref{lowhtimpliescr}, $G_{_v}^{^0} \subset G_{_v}$ is a normal subgroup of finite index.  Thus, the finite generation of $k[S]^{^{G_{_v}}}$ is reduced to checking finite generation of $k[S]^{^{G_{_v}^{^0}}}$. We may therefore replace $G_{_v}$ by $G_{_v}^{^0}$ in the proof. Further, since $G.v$ is a closed orbit, $G_{_{v,red}}$ is reductive. Therefore, by  \ref{delignestr3}, we have a surjection $M \times G_{_{v,red}} \to G_{_v}$.  Since $M$ in $G_{_v}$ is  central $k[S]^{^{G_{_v}}} = \left(k[S]^{^{{G_{_{v,red}}}}}\right)^{{M}}.$ 

Again, since $G_{_{v,red}}$ is reductive, the ring of invariants $k[S]^{^{{G_{_{v,red}}}}}$ is a finitely generated $k$-algebra. Since $M$ is a group of multiplicative type over an algebraically closed field, it is a product of $\bg_{_m}$'s and a finite group scheme. Therefore, by \cite[Theorem 1.1]{git} and \cite[Page 113]{abelian},  the ring of invariants $\left(k[S]^{^{{G_{_{v,red}}}}}\right)^{{M}}$ is also finitely generated. This proves the finite generation of $k[S]^{^{G_{_v}}}.$

The commutativity of the diagram \eqref{cartesian} now follows by the property of categorical quotients by the group $G$.

Now we check the conditions of Lemma \ref{fundlem}. We consider the given closed $G$-orbit  $F = G \times^{^{G_{_v}}} \{v\} \subset G \times^{^{G_{_v}}} S$. Then we need to check that $\phi(F)$ is a closed orbit in $V$. In fact, $\phi(F)$ is precisely the closed orbit $G.\phi(v) = G.v \subset V$. Thus we have verified conditions (2) and (3) of Lemma \ref{fundlem}.  Condition (1) is precisely the content of  \ref{lowhtimpliescr} and (4) holds since $G \times^{^{G_{_v}}} S$ has been seen to be smooth.

The isomorphism $(\text{\cursive s}, p_{_U})$ shows that the diagram \eqref{cartesian1} is cartesian. 
This completes the proof of  \ref{luna}.
\epr

\bcor\label{5.6} Let $X$ be an affine $G$-scheme  embeddable as a closed  $G$-subscheme in low height $G$-module. Let $x \in X$ be such that the orbit $G.x \subset X$ is closed. Then there exists a locally closed $G_{_x}$-invariant subscheme (a ``slice") $X_{_1}$, of $X$ with  $x \in X_{_1} \subset X$ such that the conclusions of  \ref{luna} hold for the $G$-morphism $G \times ^{^{G_{_x}}} X_{_1}  \to X $.
\ecor

\bpr This follows from \ref{luna} (which corresponds to \cite[Proposition 7.4]{bardsley}) exactly as in \cite[Proposition 7.6]{bardsley}.
\epr

The following consequence of the slice theorem has many applications so we state it here without proof. 

\bcor\label{bogocors} {\rm (see \cite[Proposition 8.5, page 312]{bardsley})} Let $F$ be an affine $G$--subvariety of ${\mathbb P}(V)$, with $V$ as a $G$--module with low height, and suppose that $F$ contains a unique closed orbit $F^{cl}$.  Then there exists a $G$-retract
$F \longrightarrow F^{cl}$.
\ecor

Let $V$ be a finite dimensional $G$--module and let 
\beqa\label{psibar}
ht_{_G}(\wedge(V)) := \mathrm{max}_{_i}\{ht_{_G}(\wedge^{^i}(V))\} 
\eeqa
We then have the following application to the theory of semistable principal bundles in positive characteristics. 

\bth\label{polystable} Let $E$ be a stable $G$-bundle  with $G$ semi-simple and $\rho_{_V}: G \rightarrow SL(V)$  be a representation such  that $p > ht_{_G}(\wedge(V))$. Then the associated bundle $E(V)$ is polystable of degree $0$. \eeth

\bpr The proof follows \cite[Theorem 9.11]{bapa2} verbatim where the notion of ``separable index" is replaced by the Dynkin height, the key ingredient being Corollary \ref{bogocors}. \epr

\bigskip

\ack \rm The first author thanks Gopal Prasad for his detailed comments on an earlier manuscript. He also thanks Brian Conrad for helpful discussions and sharing many insights and Michel Brion for his detailed comments and suggestions on an earlier manuscript. We thank the referee for his careful reading of the paper and for his comments and questions.

\bigskip

\section{Added in proof}\label{aip} In Section \ref{springerstuff}, after making suitable assumptions on $G$, we refer to Serre in \eqref{springexp} for the construction of the exponential map ${\text {exp}}:{\mathfrak g}_{_{nilp}} \stackrel{{\sim}}\longrightarrow G^{^u}$. The proof that Serre gives in (\cite[Theorem 3, page 21]{mour} is only a sketch. At his suggestion, we complete it. We work in principle over an algebraically closed field but phrase the argument so that it remains valid over more general bases.

Let $T$ be a maximal torus, $B$ a Borel subgroup containing $T$ and $U$ its unipotent radical. We assume $p \geq {\tt h}$. By this assumption on the characteristic, the Campbell-Hausdorff group law $\circ$ on $\text{Lie}(U)$ is well-defined.  The first step is the construction of a $B$-equivariant isomorphism 
\beqa\label{aip 1}
\exp : ({\text{Lie}}(U), \circ) \stackrel{{\sim}}\longrightarrow U
\eeqa
For this step we give complete proofs. The second step is to glue these isomorphisms into 
\beqa\label{aip 2}
\exp :{\mathfrak g}_{_{nilp}} \stackrel{{\sim}}\longrightarrow G^{^u}.
\eeqa
This requires the additional assumption that the simply connected covering $\widetilde{G}$ of the derived group $G'$ is \'etale over $G'$. Here we will mainly give references.
\bsem\label{aip 1.1} Let  $R^{^{+}}$ be the set of positive roots: the set of roots of $T$ acting on ${\text{Lie}}(U)$. For $\alpha$ in $R^{^{+}}$, let $U_{_\alpha}$ be the corresponding root subgroup, and let $x_{_\alpha}$ be an isomorphism from $U_{_\alpha}$ to $\bg_{_a}$. 

Fix a total order on $R^{^{+}}$. The corresponding product map from $\prod U_{_\alpha}$ to $U$ is an isomorphism of schemes: it is a morphism between smooth connected schemes which is  \'etale at the origin, bijective by \cite[13-05]{Bible}, and one  applies Zariski's main theorem as in \cite[XXII, 4.1]{sga3} (or see \cite[IV, 14.5]{Bor2}).

This isomorphism identifies $U$ with the affine space with coordinates $(x_{_\alpha})$ ($\alpha \in  R^{^{+}}$). In these coordinates, the action of $T$ on $U$ by conjugation is given by 
\begin{center}
$t$ in $T$ acts by $(x_{_\alpha}) \mapsto (\alpha(t) x_{_\alpha})$
\end{center}
This action respects the group law. It follows that the $\alpha$-coordinate of $x\cdot y$ is a linear combination of the $x^{^i}_{_\beta} y^{^j}_{_\gamma}$ for $i \beta + j\gamma = \alpha$.

We grade the polynomial algebra $\mathcal O(U)$ by defining degree of $x_{_\alpha}$ to be the {\em height} of $\alpha$ (if the expression of $\alpha$ as a linear combination of simple roots is $\sum n_{_i} \alpha_{_i}$, the height of $\alpha$ is $\sum n_{_i}$). The corresponding increasing filtration {\text {Fil}} of $\mathcal O(U)$ is invariant by $T$-conjugation as well as by the action on the right and left by $U$-translations. It is hence invariant by $B$-conjugation.

The assumption $p \geq {\tt h}$ is equivalent to $\mathcal O(U)$ being generated by $\text{Fil}_{_{p-1}}$. The definition of the exponential $\eqref{aip 1}$ is a special case of the following construction.
\esem

\bsem\label{aip 1.2} Let $U$ be an algebraic group and $\text{Fil}$ be an increasing filtration of its affine algebra. Assume that $\mathcal O(U)$ is a polynomial algebra $k[(x_{_\alpha})_{_{\alpha \in A}}]$ and that for some integers $d(\alpha) > 0$, $\text{Fil}$ is the increasing filtration deduced from the grading for which $x_{_\alpha}$ is homogeneous of degree $d(\alpha)$. Assume that the $\alpha$-coordinate of $x \cdot y$ is a polynomial of degree $\leq d(\alpha)$ in the (weighted) coordinate of $x$ and $y$. This is equivalent to the compatibility:
\beqa
\Delta(\text{Fil}_{_i}) \subset \sum_{_{j + k = i}} \text{Fil}_{_j} \otimes \text{Fil}_{_k}
\eeqa
of the filtration $\text{Fil}$ with the coproduct.

We will view $\text{Lie}(U)$ as the Lie algebra of the left invariant vector fields on $U$. Let $D_{_\alpha}$ be the element of $\text{Lie}(U)$ which is $\delta_{_\alpha}$ at $e$. If $[\epsilon]_{_\alpha}$ is the $k[\epsilon]/(\epsilon^{^2})$ valued point of $U$ with coordinate $x_{_\alpha} = \epsilon$, $x_{_\beta} = 0$ for $\beta \neq \alpha$, $D_{_\alpha}f$ is the coefficient of $\epsilon$ in $f(x.[\epsilon]_{_\alpha})$. From this it follows that $D_{_\alpha}(x_{_\beta})$ is of degree $\leq d(\beta) - d(\alpha)$, and as a consequence 
\beqa
D_{_\alpha}~ \text{Fil}_{_i} \subset \text{Fil}_{_{i - d(\alpha)}}
\eeqa
i.e. $D_{_\alpha}$ is of (descending) filtration $d(\alpha)$.

We now assume that $\mathcal O(U)$ is generated by $\text{Fil}_{_{p-1}}$, i.e. that all $d(\alpha)$ are less than $p$. Then any Lie polynomial of degree $\geq p$ in the derivations $D_{_i}$ in $\text{Lie}(U)$ vanishes on the $x_{_\alpha}$ and hence vanishes identically, implying that the Campbell-Hausdorff group law  $\circ$ on $\text{Lie}(U)$ is defined.

\bdefe\label{aip def} The exponential map ${\text {exp}}$ from $\text{Lie}(U)$ to $U$ is the map having as coordinates the truncated exponential
\beqa
{\text {exp}}(D)_{_\alpha} := \sum_{_{n < p}} \biggl(\frac{D^{^n}}{n!} (x_{_\alpha})~at~e\biggr)
\eeqa
\edefe
\blem\label{aip lemma 1} For $f$ in $\text{Fil}_{_{p-1}}$, one has 
\beqa\label{aip 3}
f(\exp(D)) = \sum_{_{n < p}} \biggl(\frac{D^{^n}}{n!} (f)~at~e\biggr)
\eeqa
\elem
\bpr This holds by the definition of ${\text {exp}}$ when $f$ is a coordinate function $x_{_\alpha}$, and it suffices to check that, for $g$ in $\text{Fil}_{_a}$ and $h$ in $\text{Fil}_{_b}$ with $a+b < p$, one has the following equality:
\beqa
\sum \frac{D^{^n}}{n!}(f g) = \sum \frac{D^{^n}}{n!}(f ) \sum \frac{D^{^n}}{n!}(g)
\eeqa
the sums being truncated at $n < p$. Both sides can be identified with sums $\sum \frac{D^{^n}}{n!}(f ) \frac{D^{^m}}{m!}(g)$, where on the left we have $n + m < p$ and on the right $n < p$ and $m < p$. If $n + m \geq p$, then either $n > a$ or $m > b$, in which case $D^{^n}(f) = 0$ (resp. $D^{^m}(g) = 0$). The terms missing on the left side therefore vanish.
\epr
\blem\label{aip lemma 2} The morphism $\exp$ is a group homomorphism
\beqa
\exp : (\mathrm{Lie}(U), \circ) \longrightarrow U.
\eeqa
\elem
\bpr If $f$ is in $\text{Fil}_{_{p-1}}$, so are its left-translates. As the derivations $D$ in $\text{Lie}(U)$ are left invariant, \eqref{aip 3} applied to the left-translate of $f$ gives:
\beqa\label{aip 4}
f(x\cdot\text{exp}(D)) = \sum_{_{n < p}} \biggl(\frac{D^{^n}}{n!} (f)~at~x\biggr)
\eeqa
hence
\beqa
f(\text{exp}(D')\cdot\text{exp}(D)) =  \sum_{_{n,m < p}} \biggl(\frac{D'^{^m}}{m!}\frac{D^{^n}}{n!} (f)~at~e\biggr)
\eeqa

In the completed free associative algebra $\bq\langle\langle X,Y \rangle\rangle$, if $H_{_d}(X,Y)$ is the Lie polynomial which is the degree $d$ component of the Campbell-Hausdorff group law, one has
\beqa
\text{exp}\biggl(\sum H_{_d}(X,Y) \biggr) = \text{exp}(X)\cdot\text{exp}(Y).
\eeqa
Modulo terms of degree $\geq p$, one still has
\beqa
\sum_{_{n \leq p}} \biggl(\sum_{_{d < p}} \frac{H_{_d}(X,Y)^{^n}}{n!} \biggr) \equiv \sum_{_{n < p}}  \frac{X^{^n}}{n!} \sum_{_{n < p}}  \frac{Y^{^n}}{n!}
\eeqa 
and this continues to hold in $\bz[\frac{1}{(p-1)!}]\langle\langle X,Y \rangle\rangle$. Substituting $D'$ and $D$ for $X$ and $Y$, we get
\beqa
f(\text{exp}(D')\cdot\text{exp}(D)) = f(\text{exp}(D'\circ D))
\eeqa
for $f$ in $\text{Fil}_{_{p-1}}$. Hence $\text{exp}(D')\cdot\text{exp}(D) = \text{exp}(D'\circ D)$.

Define $\Delta_{_a}f$ to be $f(xa) - f(x)$. For $f$ in $\text{Fil}_{_{p-1}}$, \eqref{aip 4} reads 
\beqa\label{aip 5}
\Delta_{_{\text{exp}(D)}}f = \sum_{_{0<n < p}} \biggl(\frac{D^{^n}}{n!} (f)\biggr) 
\eeqa
In $\bq[[X]]$, one has $\text{log}~\text{exp}(X) = X$. Modulo monomials of degree $\geq p$, one still has
\beqa\label{aip 6}
\sum_{_{0<n < p}} (-1)^{^n} \biggl(\sum_{_{0<m < p}} \frac{X^{^m}}{m!}  \biggr)^{^n} \big/ n \equiv X
\eeqa
and this continues to hold in $\bz[\frac{1}{(p-1)!}][[X]]$.

Substituting $D$ for $X$, in \eqref{aip 6}, we obtain  from \eqref{aip 5} that for $f$ in $\text{Fil}_{_{p-1}}$, one has
\beqa
Df = \sum_{_{0<n < p}} (-1)^{^n} \Delta^{^n}_{_{\text{exp}(D)}}f \big/ n.
\eeqa 
This gives a formula for a left inverse  $a \mapsto \text{log}~a$ of $\text{exp}:\text{Lie}(U) \longrightarrow U$; viewed as a derivation,
\beqa
\text{log}~a = \sum_{_{n < p}} (-1)^{^n} \Delta^{^n}_{_{a}} \big /n
\eeqa
when applied to $f$ in $\text{Fil}_{_{p-1}}$. This shows that the exponential is injective and hence, being an \'etale homomorphism of connected groups, it is an isomorphism.

In the case of the unipotent radical $U$ of a Borel subgroup, the given construction of the exponential map, depending only on the filtration $\text{Fil}$ of $\mathcal O(U)$, is invariant by $B$-conjugation.

\epr

\brem\label{seitzetc} (i) Let $I \subset R^{^+}$ be the set of simple roots, and for $J \subset I$, let $P_{_J}$ be the corresponding parabolic  subgroup and $U_{_J}$ its unipotent radical. Let $R^{^+}_{_J} \subset R^{^+}$ be the set of $\alpha$ in $R^{^+}$ which are linear combinations of simple roots in $J$. For $\alpha = \sum n_{_i} \alpha_{_i}$ in $R^{^+} - R^{^+}_{_J}$, define $d_{_J}(\alpha) = \sum_{_{i \notin J}} n_{_i}$. Then if all $d_{_J}(\alpha)$ are $< p$, the construction we have given of the exponential map gives a similar map
\beqa
\text{exp}: \text{Lie}(U_{_J}) \longrightarrow U_{_J}.
\eeqa
It suffices to replace $d$ by $d_{_J}$ in the earlier proofs.

(ii) More generally, let $\Gamma \subset R^{^+}$ be a closed set of positive roots, and let $U_{_\Gamma}$ be the corresponding subgroup of $U$. For any fixed total order on $\Gamma$, the corresponding product map
\beqa
\prod_{_{\alpha \in \Gamma}} U_{_\alpha}  \longrightarrow U_{_\Gamma}
\eeqa
is an isomorphism of schemes. For $\alpha$ in $\Gamma$, define $d_{_\Gamma}(\alpha)$ to be the largest number of elements of $\Gamma$  of which $\alpha$ is the sum (with repetitions allowed).

In \cite{seitz}, Seitz states that when all $d_{_\Gamma}(\alpha)$ are $< p$, an exponential map from $\text{Lie}(U_{_\Gamma})$ to $U_{_\Gamma}$ is defined. He writes (before \cite[Proposition 5.1, page 483]{seitz}) that the argument in Section 4 of \cite{serre2} yields this result. The reference is mistaken and should {\bf be replaced} by Section 2.1 of \cite{ser}.
The argument given in the preceding paragraphs, with $d$ replaced by $d_{_\Gamma}$, also yields the result. 
\erem

\esem
\bsem Let $\mathcal B$ be the variety of Borel subgroups of $G$. We have fiber spaces $X$ and $Y$ over $\mathcal B$ where the fibers at a Borel subgroup $B$ in $\mathcal B$ are respectively the unipotent radical $U$ and its Lie algebra. The given construction of the exponential map works over any base on which the primes $\ell < h$ are invertible. Over $\mathcal B$, we therefore obtain an isomorphism $\text{exp}:Y \stackrel{\sim}\longrightarrow X.$ The natural map from $X$ to $G$ (resp. $Y$ to $\text{Lie}(G)$) is proper, being the restriction to a closed subset of $\mathcal B \times G$ (resp. $\mathcal B \times \text{Lie}(G)$) of the second projection.

The reduced subscheme ${\mathfrak g}_{_{nilp}}$ (resp. $G^{^u}$)  of   ${\text{Lie}}(G)$ (resp. $G$) with points the nilpotent (resp. unipotent) elements, is the schematic image  of $Y$ in $\text{Lie}(G)$ (resp. of $X$ in $G$). We want to complete  the following $G$-equivariant diagram in the solid arrows by an isomorphism at the dotted arrow:
\beqa\label{aip 7}
\xymatrix{
Y  \ar[r]_{\sim}^{\text{exp}} \ar[d] &
X \ar[d] \\
{\mathfrak g}_{_{nilp}}  \ar@{-->}[r]&  G^{^u}
}
\eeqa
Replacing $G$ by its derived group $G'$ does not change the diagram, nor does replacing a semi-simple $G$ by its simply connected covering $\widetilde{G}$, provided $\widetilde{G}$ is \'etale over $G$. Under our assumptions this reduces us to the case where $G$ is semi-simple and simply connected. Here are references for the  existence of the dotted arrow isomorphism in \eqref{aip 7}.

(A) By \cite[6.11, 8.1]{steinberg} $G^{^u}$ is normal with a smooth dense orbit whose complement  is of codimension $\geq 2$. This holds in any characteristic for any simply connected $G$. Normality results from $G^{^u}$ being the complete intersection in $G$ defined by the equations $\chi(g) = \chi(e)$ for $\chi$ the characters of the fundamental representations, and being smooth outside codimension $2$.

(B) The same holds for ${\mathfrak g}_{_{nilp}}$ when $G$ is simply connected and $p$ good, because ${\mathfrak g}_{_{nilp}}$ is then equivariantly isomorphic to $G^{^u}$. A nice proof is in \cite[Corollary 9.3.4]{bardsley}. The isomorphism uses a representation $(V,\rho)$ such that  the pairing $\mathrm{Tr}(\rho(X)\rho(Y))$ on $\text{Lie}(G)$ is non-degenerate. It is then defined as induced by the map $g \mapsto \rho(g) - 1$ from $G$ to $\text{End}(V)$,
composed with the orthogonal projection of $\text{End}(V)$ (endowed with the symmetric bilinear form $\mathrm{Tr}(XY)$) to $\text{Lie}(G) \subset \text{End}{V}$. This does not work for $\text{SL}(n)$ for which one use $X \mapsto 1+X$.

(C) Over the dense orbits of ${\mathfrak g}_{_{nilp}}$ and $G^{^u}$, \eqref{aip 7} is a system of isomorphisms, and the composite isomorphism extends across codimension $2$ by normality.

\esem

\bigskip

\appendix{}{Coxeter number and root systems, by Zhiwei Yun}

\bprop Let $R$ be an irreducible root system associated to a simple group $G$ with Coxeter number ${\tt h}$ and let $X$ be the lattice spanned by $R$. Let $\phi:X \to \br/\bz$ be a homomorphism. Then there exists a basis $B$ for $R$ such that if $\alpha \in R$ satisfies $\phi(\alpha) \in \big(0, 1/{\tt h}\big) ~mod~ \bz$, then $\alpha$ is positive with respect to $B$. \eprop

\bpr Pick any basis $B_{_0}  = \{\alpha_{_1}, \ldots , \alpha_{_r}\}$ of $R$. We lift $\phi$ to a linear map $y:X \to \br$ which is viewed as a point in $V = {\text {Hom}}(X, \br)$. 

Let $W$ be the Weyl group of $R$. Then there is an action of the extended Weyl group $ {\text {Hom}}(X, \bz) \rtimes W$ on $V$ and the validity of the statement for the linear map $y$ is invariant for this action. Therefore we may assume that $y$ lies in the fundamental alcove attached to $B_{_0}$, i.e. 
\beqa
\alpha_{_j}(y) \geq 0, ~~~~\forall 1 \leq j \leq r \\
\theta(y) \leq 1
\eeqa
where $\theta$ is the highest root with respect to $B_{_0}$. Then for all $\alpha \in R$, we have 
\beqa
\mid \alpha(y) \mid \leq 1.
\eeqa
We {\em claim} that there exists a basis $B$ of $R$ such that for any $\beta \in R$ with 
\beqa\label{one}
\beta(y) \in \big(-1, -1 + 1/{\tt h} \big) \cup \big(0, 1/{\tt h}\big),
\eeqa
$\beta$ positive with respect to $B$.

Write $\theta = \sum_{_i} n_{_i} \alpha_{_i}$ and let $\alpha_{_0} = 1 - \theta$ and $n_{_0} = 1$ as usual. Then we have the equality $\sum_{_{i = 0}}^{^r} n_{_i} \alpha_{_i} = 1$, viewed as affine functions on $V$ and we have 
\beqa
{\tt h} = \sum_{_{i = 0}}^{^r} n_{_i}.
\eeqa
Therefore, for some $0 \leq i \leq r$, we have $\alpha_{_i}(y) \geq 1/{\tt h}$.

If $i = 0$, then $\theta(y) \leq 1 - 1/{\tt h}$ and in this case, we can take $B = B_{_0}$. In fact, for any negative root $-\beta$ with respect to $B_{_0}$, we have $0 \geq -\beta(y) \geq 1 - 1/{\tt h}$, so that $-\beta$ cannot satisfy the condition \eqref{one}.

If $i > 0$, let $v_{_i}$ be the vertex of the fundamental alcove opposite to the root hyperplane $\alpha_{_i}$. So $\alpha_{_j}(v_{_i}) = 0$ if $j \neq i$ and $\alpha_{_i}(v_{_i}) = 1/n_{_i}$. We take a basis $B$ of $R$ such that the vector $y - v_{_i}$ is dominant, i.e. for any $\beta \in R$ positive for $B$, we have the relation: $\beta(y) \geq \beta(v_{_i})$.

Such a basis $B$ of $R$ satisfies the requirements of the claim. In fact, let $\beta \in R$ be negative with respect to $B$.

Write $\beta = \sum_{_{j = 1}}^{^r} m_{_j} \alpha_{_j}$ for some $-n_{_j} \leq m_{_j} \leq n_{_j}$. Then, $\beta(v_{_i}) = m_{_i}/n_{_i}$. On the other hand, $\beta(y) = \sum_{_{j = 1}}^{^r} m_{_j} \alpha_{_j}(y)$. Since $\beta$ is negative for $B$, we have 
\beqa\label{two}
\beta(y) \leq \beta(v_{_i}) = m_{_i}/n_{_i}
\eeqa
and we have two cases:
\begin{itemize}
\item If $\beta > 0$, with respect to $B_{_0}$, then $\beta(y) \geq 0$ and $0 \leq m_{_j} \leq n_{_j}$ for all $j$. If $m_{_i} = 0$, then  \eqref{two} implies that $\beta(y) = 0$ and so $\beta$ fails to satisfy \eqref{one}. If $m_{_i} \geq 1$, then $\beta(y) \geq \alpha_{_i}(y) \geq 1/{\tt h}$, again not satisfying \eqref{one}.
\item If $\beta < 0$, with respect to $B_{_0}$, then $-1 \leq \beta(y) \leq 0$ and $-n_{_j} \leq m_{_j} \leq 0$ for all $j$. If $m_{_i} = -n_{_i}$, then $\beta(v_{_i}) = -1$ which forces $\beta(y) = -1$ by \eqref{two}. If $m_{_i} \geq -n_{_i} + 1$, then comparing $\mid \beta(y) \mid$ with $\mid \theta(y) \mid$, it misses at least one copy of $\alpha_{_i}(y)$, therefore $\mid \beta(y) \mid \leq \theta(y) - \alpha_{_i}(y) \leq 1 - 1/{\tt h}$ and yet again, $\beta(y)$ does not satisfy \eqref{one}.
\end{itemize}

This proves the claim and the proposition.
\epr

\bibliographystyle{amsalpha}
\bibliographymark{References}

\end{document}